\documentclass[num-refs]{BiB_Styles/wiley-article}

\usepackage{lineno}
\usepackage{amsmath}
\usepackage{amssymb}
\usepackage{subfig}
\usepackage{graphicx}
\usepackage{stmaryrd}
\usepackage{tikz}
\usepackage{natbib}

\newcommand{\llangle}{\left\langle}
\newcommand{\rrangle}{\right\rangle}

\newcommand{\wb}{\mathbb{W}}
\newcommand{\detj}{\det\mathbf{J}}

\newcommand{\uadvc}{\overline{\widehat{\mathbf{u}}}^{1/2}}
\newcommand{\uadv}{\overline{\mathbf{u}}^{1/2}}

\DeclareFontFamily{OMX}{MnSymbolE}{}
\DeclareSymbolFont{MnLargeSymbols}{OMX}{MnSymbolE}{m}{n}
\SetSymbolFont{MnLargeSymbols}{bold}{OMX}{MnSymbolE}{b}{n}
\DeclareFontShape{OMX}{MnSymbolE}{m}{n}{
    <-6>  MnSymbolE5
   <6-7>  MnSymbolE6
   <7-8>  MnSymbolE7
   <8-9>  MnSymbolE8
   <9-10> MnSymbolE9
  <10-12> MnSymbolE10
  <12->   MnSymbolE12
}{}
\DeclareFontShape{OMX}{MnSymbolE}{b}{n}{
    <-6>  MnSymbolE-Bold5
   <6-7>  MnSymbolE-Bold6
   <7-8>  MnSymbolE-Bold7
   <8-9>  MnSymbolE-Bold8
   <9-10> MnSymbolE-Bold9
  <10-12> MnSymbolE-Bold10
  <12->   MnSymbolE-Bold12
}{}

\DeclareMathDelimiter{\lllangle}{\mathopen}%
                     {MnLargeSymbols}{'164}{MnLargeSymbols}{'164}
\DeclareMathDelimiter{\rrrangle}{\mathclose}%
                     {MnLargeSymbols}{'171}{MnLargeSymbols}{'171}

\makeatletter
\def\blfootnotetext{\gdef\@thefnmark{}\@footnotetext}
\makeatother

\papertype{Original Article}

\title{A mixed finite-element, finite-volume, semi-implicit discretisation
for atmospheric dynamics: Spherical geometry}

\author[1]{Thomas~Melvin}
\author[1]{Ben~Shipway}
\author[1]{Nigel~Wood}
\author[1$\dagger$]{Tommaso~Benacchio} 
\author[1]{Thomas~Bendall}
\author[1]{Ian Boutle}
\author[1]{Alex~Brown}
\author[1]{Christine Johnson}
\author[1]{James~Kent}
\author[1]{Stephen~Pring}
\author[1]{Chris Smith}
\author[1]{Mohamed~Zerroukat}
\author[2]{Colin~Cotter}
\author[3]{John~Thuburn}
\affil[1]{Met Office, Exeter, United Kingdom}
\affil[2]{Imperial College, London, United Kingdom}
\affil[3]{University of Exeter, United Kingdom}
\affil[$\dagger$]{Current Affiliation: Danish Meteorological Institute, Copenhagen, Denmark}

\corraddress{Thomas Melvin, Met Office, FitzRoy Road, Exeter EX1 3PB, United Kingdom}
\corremail{Thomas.Melvin@metoffice.gov.uk}

\fundinginfo{NERC through grants NE/K006762/1 \&
NE/K006789/1. \\
EPSRC through grants EP/L016613/1 \& EP/R029423/1.}

\blfootnotetext{\textcopyright \hspace{0.5mm} Crown copyright 2024. Reproduced with the permission of the Controller of HMSO. Published by John Wiley and Sons, Ltd.}

\runningauthor{Melvin et al.}

\begin{document}

\maketitle
\begin{abstract}
The reformulation of the Met Office's dynamical core for weather and climate prediction previously described by the authors is extended to spherical domains using a cubed-sphere mesh. This paper updates the semi-implicit mixed finite-element formulation to be suitable for spherical domains. In particular the finite-volume transport scheme is extended to take account of non-uniform, non-orthogonal meshes and uses an advective-then-flux formulation so that increment from the transport scheme is linear in the divergence. The resulting model is then applied to a standard set of dry dynamical core tests and compared to the existing semi-implicit semi-Lagrangian dynamical core currently used in the Met Office's operational model.  
\end{abstract}
\section{Introduction}

\label{sec:Introduction}

At the centre of all weather and climate models lies the dynamical core. The dynamical core approximates the fluid dynamical motion that is resolved by the model mesh and is  coupled to models for unresolved processes such as the boundary layer and non-fluid processes such as radiation. The dynamical core is required to be accurate, stable and efficient for the scales of motion that it simulates. Fundamental to achieving these properties is the choice of model mesh. 
This choice can result in a number of features that need to be addressed by the numerical scheme, such as resolution clustering, non-orthogonality, grid imprinting and computational modes; see \citet{Staniforth12} for a more detailed discussion.

Modern supercomputers consist of a greatly increasing number of (increasingly heterogeneous) processors and in order to take advantage of this computational resource the dynamical core needs to make efficient use of memory management and communication processes (\citet{Lawrence18}). This has led to a shift away from the regular Latitude-Longitude mesh (which, due to convergence of the meridians at the poles, leads to computational bottlenecks) and towards some form of quasi-uniform horizontal mesh. \citet{Staniforth12} detailed a number of desirable properties that any numerical scheme designed for dynamical cores should exhibit, including mimetic/compatible properties, at least 2nd order accuracy and minimal grid imprinting.  Achieving these properties is non-trivial, particularly on a non-orthogonal meshes, such as is typically the case for quasi-uniform meshes on the sphere. 

\citet{cottershipton12, cotterthuburn14} and \citet{thuburncotter15} developed a family of compatible mixed finite-element methods for the shallow-water equations on the sphere where orthogonality of the underlying mesh is not required to achieve good accuracy and hence they are well suited to a quasi-uniform mesh. This family of schemes was applied to a variety of icosahedral and cubed-sphere meshes. The mixed finite-element approach was extended into three-dimensions by \citet{natale2016} and \citet{melvin19} who presented an application of this mixed finite-element model in Cartesian geometry using hexahedral elements and coupled to a finite-volume transport scheme and semi-implicit timestepping. \citet{kent22} then extended the formulation of \citet{melvin19} to the shallow water equations on a cubed-sphere horizontal mesh. See \citet{cotter2023compatible} for a survey of mixed finite-element methods for geophysical modelling. 

The compatible mixed finite-element approach fulfils many of the
 desirable properties detailed in \citet{Staniforth12} for the design of a dynamical core. Importantly it provides discrete analogues of certain continuous vector calculus identities (such as $\nabla\times\nabla \psi\equiv 0$ and $\nabla\cdot\nabla\times\mathbf{v}\equiv 0$ for all scalar $\psi$ and vector $\mathbf{v}$) as well as sharing many of the good wave dispersion properties of the widely used Arakawa C-grid staggering (\citet{arakawa77}). This scheme meets the necessary conditions for the absence of computational modes (such as a 2:1 ratio of horizontal velocity degrees of freedom to pressure degrees of freedom). \citet{natale2016} \& \citet{melvin18} showed how to create function spaces that mimic the Charney-Phillips grid staggering in the vertical direction which is desirable due to the absence of computational modes and good wave dispersion properties, (\citet{thuburnwoolings05}).

In this paper the formulation used by \citet{melvin19} and \citet{kent22} is extended to three-dimensional spherical domains on a cubed-sphere mesh. Principle differences from these models are given in Section \ref{sec:formulation} and the governing equations that are used are revisited in Section \ref{sec:Continuous-equations}. The wave dynamic components of the model are spatially discretised using the mixed finite-element method of \citet{cottershipton12} and the temporal discretisation uses an iterated-semi-implicit scheme inspired by that of \citet{wood14}, seeking to maintain the temporal accuracy and long timestep stability of that model, Sections \ref{sec:SI-temporal-discretisation}-\ref{sec:Spatial-discretization}. Mappings from the computational space to the Equiangluar cubed-sphere mesh used by the finite-element scheme in this study are given in Section \ref{sec:mesh}. As in \citet{melvin19} and \citet{kent22} the finite-element wave dynamics model is coupled to an explicit finite-volume scheme for the transport terms, described in Section \ref{sec:transport}, which is applied to all model variables. Compared to the scheme of \citet{melvin19} and \citet{kent22}, a directionally split method-of-lines scheme is used here to improve the efficiency of the transport scheme, and an advective-then-flux approach is used to improve stability of conservative transport in this context. Following \citet{wood14} the iterative timestep is split into outer (transport) and inner (nonlinear) loops and at each iteration a linear system inspired by a semi-implicit formulation is solved as in \citet{maynard20}, outlined in Section \ref{sec:timestepping}. To assess the model's behaviour in a range of flow regimes it is applied to a number of dynamical core tests from the literature; the results are presented in Section \ref{sec:Results} and compared to the semi-implicit semi-Lagrangian model of \citet{wood14}. Finally, conclusions are summarised in Section \ref{sec:Discussion}.

\section{Model Formulation}
\label{sec:formulation}
The model formulation closely follows that of \citet{melvin19} and the formulation is revisited in the following sections. However there a number of differences that are highlighted here and discussed, together with their motivation, in more detail later:
\begin{enumerate}
    \item Changes to the finite-volume transport scheme:
      \begin{itemize}        
        \item The momentum equation is reformulated in the advective form instead of the vector invariant form and the advection terms are handled by the explicit finite-volume advection scheme instead of the semi-implicit mixed finite-element method used by \citet{melvin19}  (Section \ref{sec:Spatial-discretization})
        \item The finite-volume transport scheme is temporally split between vertical and horizontal directions using a Strang splitting method (Section \ref{sec:transport_strang});
        \item The conservative transport scheme uses an advective-then-flux formulation (Section \ref{sec:adv-then-flux}).
        \item The polynomial reconstruction used in the method-of-lines advection scheme uses a two-dimensional horizontal (\citet{kent22}) and one-dimensional vertical reconstruction (Section \ref{sec:spatial-reconstruction});
      \end{itemize}
    \item The geopotential is placed in the $\wb_3$ space instead of $\mathbb{W}_0$ (Section \ref{sec:Spatial-discretization});
    \item The Jacobian $\mathbf{J}$ mapping from the computational space to the physical space is computed with a semi-analytic expression via an intermediate spherical coordinate system (Section \ref{sec:mesh});
    \item The equation of state is sampled at nodal points of the $\wb_3$ degrees of freedom instead of being solved in the weak form. (Section \ref{sec:discrete-linear-system})
\end{enumerate}

\section{Continuous equations}

\label{sec:Continuous-equations}

The Euler equations for a dry perfect gas in a rotating frame are
\begin{eqnarray}
\frac{\partial\mathbf{u}}{\partial t} & = &  -\left(\mathbf{u}\cdot\nabla\right)\mathbf{u} + \mathbf{S},\label{eq:velocity_cont}\\
\frac{\partial\rho}{\partial t} & = & -\nabla\cdot\left(\rho\mathbf{u}\right),\label{eq:continuity_cont}\\
\frac{\partial\theta}{\partial t} & = & -\mathbf{u}\cdot\nabla\theta,\label{eq:thermo_cont}
\end{eqnarray}
where $\mathbf{S}\equiv -2\boldsymbol{\Omega}\times\mathbf{u}-\nabla\Phi-c_{p}\theta\nabla\Pi$, together with the nonlinear equation of state \begin{equation}
\Pi^{\left(\frac{1-\kappa}{\kappa}\right)}=\frac{R}{p_{0}}\rho\theta.\label{eq:eos_cont}\end{equation}
The velocity vector is $\mathbf{u}$; $\boldsymbol{\Omega}$ is the rotation vector;
$\Phi$ is the geopotential; $c_{p}$
is the specific heat at constant pressure; $\theta$ is the potential
temperature, related to temperature through $T=\theta\Pi$; $\Pi=\left(p/p_{0}\right)^{\kappa}$
is the Exner pressure with $p$ pressure and $p_{0}$ a constant reference
pressure; $R$ is the gas constant per unit mass; $\kappa\equiv R/c_{p}$;
and $\rho$ is the density.

These equations are solved on a spherical shell subject to the boundary condition of zero mass flux
through the top and bottom boundaries of the domain.

\section{Overview of the spatio-temporal discretisation }

\label{sec:SI-temporal-discretisation}

The temporal discretisation of the equations follows that of \citet{melvin19} and is inspired by an iterative-semi-implicit
semi-Lagrangian discretisation such as that used in \citet{wood14}.
In this scheme the advective terms (including, in contrast to \citet{melvin19}, the advection terms in the momentum equation) are handled using an explicit Eulerian scheme which acts on an intermediate update of the wave dynamics terms (see Section \ref{sec:transport_state}) for the variables (instead of the time level $n$ terms in \citet{melvin19}), either in flux form for the continuity equation or advective form for potential temperature and momentum equations. All
other terms are handled using an iterative-implicit temporal discretisation. The momentum equation is recast from the vector invariant form of \citet{melvin19} to the advective form, resulting in a consistent discretisation of all transport terms using the finite-volume transport scheme of Section \ref{sec:transport} due to the explicit presence of a transport term for the wind field.

With the addition of an implicit damping layer applied to the vertical component of the velocity vector \eqref{eq:velocity_cont}-\eqref{eq:eos_cont} are discretised to give
\begin{eqnarray}
\delta_{t}\mathbf{u} &=& -\overline{\mu\left(\frac{\mathbf{u}\cdot\mathbf{n}}{\mathbf{z}\cdot\mathbf{n}}\right)\mathbf{z}}^1-\boldsymbol{\mathcal{A}}\left(\mathbf{u}^{p},\uadv\right)  + \overline{\mathbf{S}}^{\alpha},\label{eq:velocity_disc_temp}\\
\delta_{t}\rho &=&-\nabla\cdot\boldsymbol{\mathcal{F}}\left(\rho^{p},\uadv\right),\label{eq:continuity_disc_temp}\\
\delta_{t}\theta &=& -\mathcal{A}\left(\theta^{p},\uadv\right),\label{eq:thermo_disc_temp}\\
\overline{\Pi^{\left(\frac{1-\kappa}{\kappa}\right)}}^1&=&\overline{\frac{R}{p_{0}}\rho\theta}^1,\label{eq:eos_disc_temp}
\end{eqnarray}
where,
for a generic scalar or vector variable $F$,\begin{equation}
\delta_{t}F\equiv\frac{F^{n+1}-F^{n}}{\Delta t},\qquad\overline{F}^{\alpha}\equiv\alpha F^{n+1}+\left(1-\alpha\right)F^{n}.\label{eq:delta_t_def}\end{equation}
The parameter $\alpha$ is a temporal off-centring parameter and the superscripts $n$ and $n+1$ indicate approximations at time $n\Delta t$ and $\left(n+1\right)\Delta t$ respectively.
The advecting velocity $\uadv$ is therefore a centred Eulerian average in time and, in contrast to \citet{melvin19}, the consistent metric modification of $\uadv$ (their Section 5.3.3) is not used here. $\boldsymbol{\mathcal{F}}\left(q^{p},\uadv\right)$
is the time-averaged flux and $\mathcal{A}\left(q^{p},\uadv\right)$
is the time-averaged advection tendency of a scalar $q^p$ and 
$\boldsymbol{\mathcal{A}}\left(\mathbf{v}^{p},\uadv\right)$
is the time-averaged advection tendency of a vector $\mathbf{v}^{p}$, where superscripts $p$ indicate an intermediate wave dynamics state of a scalar $q$ or vector field $\mathbf{v}$. See Section \ref{sec:transport_state} for more details how the transported states $q^p$ and $\mathbf{v}^{p}$ are computed.

All terms are discretised in space using the mixed finite-element
scheme described in section \ref{sec:Spatial-discretization}, except
for $\boldsymbol{\mathcal{F}}$, $\mathcal{A}$ and $\boldsymbol{\mathcal{A}}$ which are discretised
using the finite-volume scheme described in section \ref{sec:transport}.


\section{The mixed finite-element discretisation}

\label{sec:Spatial-discretization}

The finite-element spaces used in the spatial discretisation are the same as those in \citet{melvin19}.
Each variable is represented in an appropriate function space in a hexahedral element: 
\begin{itemize}
\item $\mathbf{u}\in\wb_{2}$: The Raviart-Thomas $RT_l$ space of vector functions of degree $l$ tangential to and discontinuous across an element facet and degree $l+1$ normal to and continuous across an element facet;
\item $\Pi,\,\rho,\,\Phi\in\wb_{3}$: The $Q^{DG}_{l}$ space of scalar functions built from the tensor product of degree $l$ polynomials that are discontinuous at element boundaries.
\item $\theta\in\wb_\theta$: The space of scalar functions based on the vertical part of $\mathbb{W}_{2}$
\item The components of the coordinate field $\chi_i\in\wb_\chi,\,i=1,2,3$: The $Q^{DG}_{m}$ space of scalar functions.
\end{itemize}

As in \citet{melvin19} equations \eqref{eq:velocity_disc_temp}-\eqref{eq:eos_disc_temp} are transformed from cells in the physical space to a single reference cell in the computational space. The physical space and the computational space are both assumed to be Cartesian spaces. The associated metric tensor is the identity and is therefore dropped in the following formulation. In principle it would be possible to use a spherical coordinate system as the physical coordinates, however this would introduce a non-diagonal metric tensor in the mapping from physical to computational space and would modify the transformations used below. 

Given a mapping (Piola transformation), $\boldsymbol{\phi}$, between physical space (denoted by undressed variables) $\boldsymbol{\chi}$ and computational space (denoted by dressed, $\widehat{\hphantom{v}}$ variables) $\widehat{\boldsymbol{\chi}}$ such that $\boldsymbol{\chi} = \boldsymbol{\phi}(\widehat{\boldsymbol{\chi}})$, with a Jacobian $\mathbf{J}\equiv\partial\boldsymbol{\phi}\left(\widehat{\boldsymbol{\chi}}\right)/\partial\widehat{\boldsymbol{\chi}}$, variables in each of the four spaces transform from computational space to physical space according to the following rules (see \citet{rognes10} and references within for more details)
\begin{itemize}
    \item $\wb_{2}$: $\mathbf{v}=\mathbf{J}\widehat{\mathbf{v}}/\det\mathbf{J}$;
    \item $\wb_{3}$: $\sigma=\widehat{\sigma}$ and $\nabla\cdot\mathbf{v}=\widehat{\nabla}\cdot\widehat{\mathbf{v}}/\det\mathbf{J}$ for $\mathbf{v}\in\wb_2$;
    \item $\wb_\theta$: $w = \widehat{w}$;
    \item $\wb_\chi$: $\zeta = \widehat{\zeta}$
\end{itemize}
Note that, as in \citet{melvin19}, to avoid problems with not being able to exactly integrate the weak form of divergence, the rehabilitation method of \citet{bochev_ridzal08} is used. This modifies the mapping of $\sigma\in\wb_3$ from $\sigma=\widehat{\sigma}/\detj$ to $\sigma = \widehat{\sigma}$ which results in a weak divergence $\int \widehat{\sigma}\widehat{\nabla}\cdot\widehat{\mathbf{v}}$ that can be exactly integrated for general cell shapes, see \citet{natale2016} \& \citet{melvin19} for more details. In contrast to \citet{melvin19}, the $\wb_{0}$ and $\wb_{1}$ spaces are not used in this formulation. Placing $\Phi\in\wb_3$ gives a more compact stencil for the geopotential gradient that matches that of the pressure gradient. This leads to a small improvement in the model's discrete representation of quasi-hydrostatic balance.

%

%

\subsection{Discrete equations using the computational cell}

The discretisations of the  continuity \eqref{eq:continuity_disc_temp} and thermodynamic \eqref{eq:thermo_disc_temp} equations follow \citet{melvin19} and so are not repeated here. The momentum equation \eqref{eq:velocity_disc_temp} is now cast in the advective form and therefore the discretisation of the transport terms has altered. The form used here is obtained by multiplying \eqref{eq:velocity_disc_temp} by a test function $\mathbf{v}\in\wb_2$, transforming to the computational space and integrating over the domain, giving

\begin{eqnarray}
\llangle \mathbf{J}\widehat{\mathbf{v}},\frac{\mathbf{J}\delta_{t}\widehat{\mathbf{u}}}{\det\mathbf{J}}\rrangle &=&
- \overline{\llangle \mathbf{J}\widehat{\mathbf{v}},\mu\left(\frac{\widehat{\mathbf{u}}\cdot\widehat{\mathbf{n}}_{b}}{\widehat{\mathbf{z}}_{b}\cdot\widehat{\mathbf{n}}_{b}}\right)\frac{\mathbf{J}\widehat{\mathbf{z}}_{b}}{\det\mathbf{J}}\rrangle }^{1} + R_u^A + \overline{R_u}^\alpha,\label{eq:velocity_weak_pullback}\\
R_u^A  & \equiv & -\llangle\mathbf{J}\widehat{\mathbf{v}},\boldsymbol{\mathcal{A}}\left(\widehat{\mathbf{u}}^p,\uadvc\right)\rrangle, \label{eq:rhs-u-adv}\\
\overline{R_u}^\alpha & \equiv & -\overline{\llangle \mathbf{J}\widehat{\mathbf{v}},2\boldsymbol{\Omega}\times\frac{\mathbf{J}\widehat{\mathbf{u}}}{\det\mathbf{J}}\rrangle }^{\alpha} +\overline{\left\langle \widehat{\nabla}\cdot\widehat{\mathbf{v}},\widehat{\Phi}\right\rangle }^{\alpha}\nonumber \\
  &  & -\overline{\left\lllangle\llbracket c_{p}\widehat{\theta}\widehat{\mathbf{v}}\rrbracket,\left\{ \widehat{\Pi}\right\} \right\rrrangle}^{\alpha} +\overline{\llangle c_{p}\widehat{\theta}\widehat{\nabla}_{\mathsf{C}}\cdot\widehat{\mathbf{v}}+\widehat{\mathbf{v}}\cdot\widehat{\nabla}_{\mathsf{C}}\left(c_{p}\widehat{\theta}\right),\widehat{\Pi}\rrangle }^{\alpha},
\end{eqnarray}
where $\llangle\llangle\cdot\rrangle\rrangle$ denotes the surface integrals
over the collection of all cell faces evaluated in the computational space
and $\llbracket\cdot\rrbracket$ and $\left\{\cdot\right\}$ indicate the jump in its argument and the value of its argument accross a cell face, respectively (see \citet{melvin19}, section 4.1.1. for details).

The equation of state \eqref{eq:eos_disc_temp} is sampled at nodal points of the finite-element scheme rather than solved in its weak form as in \citet{melvin19}. This is motivated by a desire for \eqref{eq:eos_disc_temp} to hold exactly as a diagnostic relationship between the Exner pressure, potential temperature and density.

\section{Mesh Mappings}
\label{sec:mesh}

The discretisation presented in the previous section is valid for a general three-dimensional hexahedral mesh
and was applied to a uniform horizontally biperiodic domain in \citet{melvin19}. Here
spherical shell domains are considered.

The computational space consists of a single unit cell with Cartesian coordinates $\widehat{\boldsymbol{\chi}}\in\left[0,1\right]^3$ and mappings $\boldsymbol{\chi}=\boldsymbol{\phi}\left(\widehat{\boldsymbol{\chi}}\right)$ are introduced to map this computational cell to each cell in the physical mesh.

The physical mesh coordinates are parametrised as a finite-element field $\boldsymbol{\chi}\in\wb_\chi$ as in \citet{melvin19} with a piecewise polynomial representation of degree $m$. However, representing a spherical manifold with a piecewise polynomial introduces discretisation errors that depend upon the degree $m$ and in order to accurately represent the surface of the sphere a high degree $\wb_\chi$ space is required (\citet{kent22}). To avoid this, an alternative approach used here is to map via an intermediate spherical coordinate system $\boldsymbol{\xi}$ such that the Jacobian is given by
\begin{equation}
    \mathbf{J} \equiv \frac{\partial\boldsymbol{\phi}\left(\widehat{\boldsymbol{\chi}}\right)}{\partial\widehat{\boldsymbol{\chi}}}\equiv \frac{\partial\boldsymbol{\chi}}{\partial\boldsymbol{\xi}}\frac{\partial\boldsymbol{\xi}}{\partial\widehat{\boldsymbol{\chi}}},\label{eq:spherical-Jacobian}
\end{equation}
where $\boldsymbol{\xi}$ can be chosen such that the transformation of a cell from computational space to physical space, $\partial\boldsymbol{\xi}/\partial\widehat{\boldsymbol{\chi}}$, can be accurately represented by $\boldsymbol{\xi}\in\wb_\chi$ with a low degree $\wb_\chi$ space. The coordinate transformation from $\boldsymbol{\xi}$ to $\boldsymbol{\chi}$,   $\partial\boldsymbol{\chi}/\partial\boldsymbol{\xi}$ is chosen such that it has a known analytic form that captures the spherical nature of the manifold. The form of $\mathbf{J}$ used in \citet{melvin19} is recovered by setting
$\boldsymbol{\xi}=\boldsymbol{\chi}$ in which case $\partial\boldsymbol{\chi}/\partial\boldsymbol{\xi}$ becomes the identity matrix.

\begin{figure}
\centering
\includegraphics[scale=0.3]{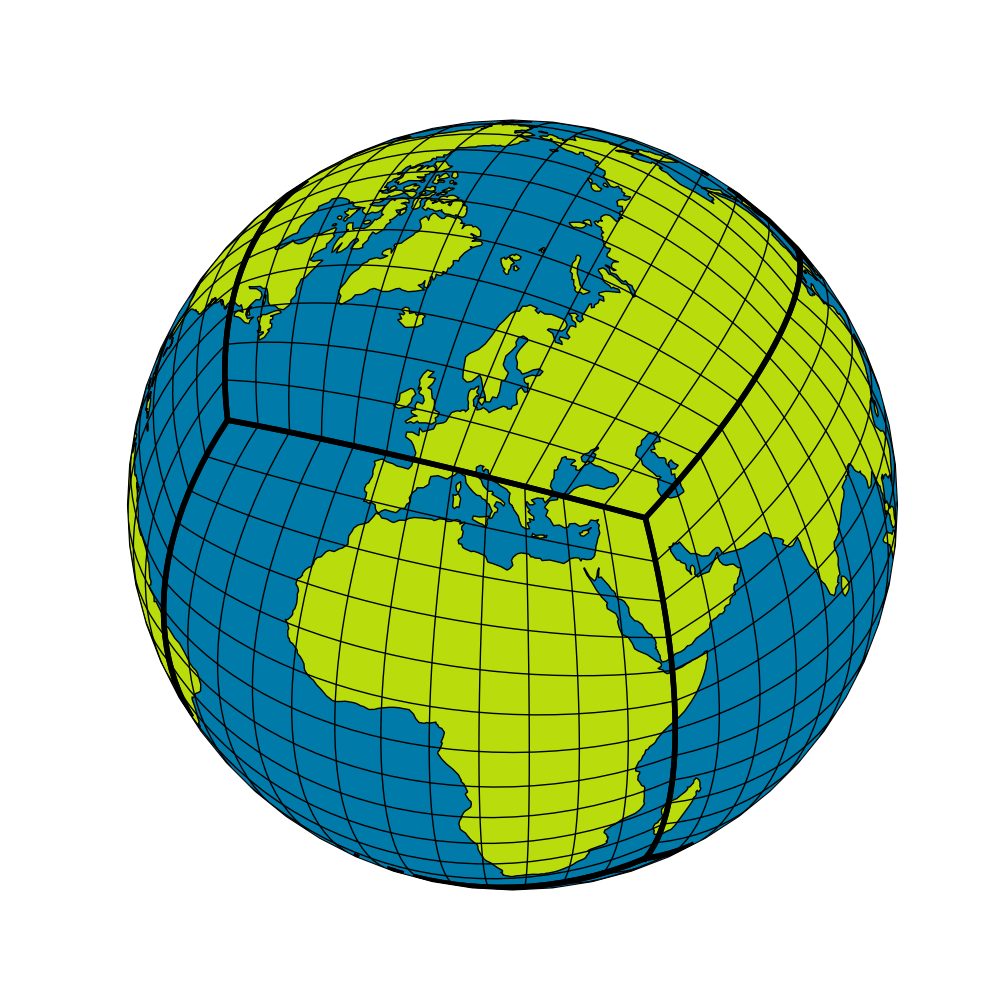}
\caption{$C12$ cubed-sphere horizontal mesh with $6x12x12$ cells.}
\label{fig:cubed-sphere}
\end{figure}

The horizontal mesh used here is an equi-angular cubed-sphere (\citet{ronchi_etal_1996}), 
Figure \ref{fig:cubed-sphere}.
Vertically the mesh is extruded as described in \citet{adams_etal19}. The mesh resolution is denoted as $CnLm$ where 
$n$ is the number of cells along one edge of a panel and the $m$ is the number of vertical layers such that 
there are $6n^2$ model columns and $6n^2m$ cells in the three-dimensional mesh.

A geocentric Cartesian coordinate system $\boldsymbol{\chi}\equiv\left(X,Y,Z\right)$ is used where $\left(X,Y,Z\right)=0$ is the centre of the sphere of radius $a=\sqrt{X^2+Y^2+Z^2}$.
Alongside the Cartesian coordinates $\boldsymbol{\chi}$ a spherical coordinate system $\boldsymbol{\xi}\equiv\left(\xi,\eta,r\right)$ is used with angular variables $\left[\xi,\eta\right]\in\left[-\pi/4,\pi/4\right]$ on each panel such that lines of constant $\xi$ and $\eta$ are angularly equidistant great circles on each panel and $r$ is the radial distance from the centre of the sphere.
As an example, with this choice of coordinates and mesh the components of the Jacobian $\mathbf{J}\equiv\frac{\partial\boldsymbol{\chi}}{\partial\boldsymbol{\xi}}\frac{\partial\boldsymbol{\xi}}{\partial\widehat{\boldsymbol{\chi}}}$ in each cell $c$ with spacing $\Delta\xi_c$ and $\Delta\eta_c$ in the $\xi$ and $\eta$ direction respectively and with a constant slope in height above the surface of the sphere of $\delta_\xi\equiv\Delta r_c/\Delta \xi_c$ in the $\xi-$direction and $\delta_\eta\equiv\Delta r_c/\Delta \eta_c$ in the $\eta-$direction are 
\begin{equation}
    \frac{\partial\boldsymbol{\xi}}{\partial\widehat{\boldsymbol{\chi}}} = \begin{pmatrix}
    \Delta\xi_c & 0 & 0 \\
    0 & \Delta\eta_c & 0 \\
    \Delta\xi_c \tan \delta_\xi & \Delta\eta_c \tan \delta_\eta & \Delta r_c
    \end{pmatrix}.\label{eq:dxidchi}
\end{equation}
In practice, to maintain flexibility of the scheme and to facilitate the inclusion of arbitrary orography, this component of the Jacobian is computed numerically with $\boldsymbol{\xi}\in\mathbb{W}_\chi$. 

The second component of the Jacobian, $\partial\boldsymbol{\chi}/\partial\boldsymbol{\xi}$, transforms the spherical coordinate $\boldsymbol{\xi}$ into the Cartesian coordinate $\boldsymbol{\chi}$. Following \citet{nair05} the basis vectors for a panel of the Equiangular cubed-sphere are
\begin{eqnarray}
    \mathbf{e}_\xi &=& (\frac{r}{\varrho^3})(1 + t_\xi^2)\left[-t_\xi,\,(1 + t_\eta^2),\,-t_\xi t_\eta\right],\\
    \mathbf{e}_\eta &=& (\frac{r}{\varrho^3})(1 + t_\eta^2)\left[-t_\eta,\,-t_\xi t_\eta,\,(1 + t_\xi^2)\right],\\
    \mathbf{e}_r &=& \frac{1}{\varrho}\left[1,\,t_\xi,\,t_\eta\right], 
\end{eqnarray}
where
\begin{equation}
    t_\xi = \tan(\xi),\quad
    t_\eta = \tan(\eta),\quad
    \varrho = \sqrt{1 + t_\xi^2 + t_\eta^2}.
\end{equation}
The second component of the Jacobian mapping is then
\begin{equation}
    \frac{\partial\boldsymbol{\chi}}{\partial\boldsymbol{\xi}} = R_i\left[\mathbf{e}^T_\xi,\,\mathbf{e}^T_\eta,\,\mathbf{e}^T_r\right],
\end{equation}
where $R_i$ is a rotation matrix for panel $i=1,\ldots,6$ of the cubed-sphere that will translate and rotate $\left[\mathbf{e}^T_\xi,\,\mathbf{e}^T_\eta,\,\mathbf{e}^T_r\right]$ such that the union of all 6 panels form a spherical shell. The rotation matrix for each panel is given by
\begin{equation}
  \begin{matrix}
  R_1 = \begin{pmatrix}1 & 0 & 0\\ 0 & 1 & 0\\ 0 & 0 & 1\end{pmatrix},&
  R_2 = \begin{pmatrix}0 & -1 & 0\\ 1 & 0 & 0\\ 0 & 0 & 1\end{pmatrix},&
  R_3 = \begin{pmatrix}-1 & 0 & 0\\ 0 & 0 & 1\\ 0 & 1 & 0\end{pmatrix},\\
  R_4 = \begin{pmatrix}0 & 0 & -1\\ -1 & 0 & 0\\ 0 & 1 & 0\end{pmatrix},&
  R_5 = \begin{pmatrix}0 & 0 & -1\\ 0 & 1 & 0\\ 1 & 0 & 0\end{pmatrix},&
  R_6 = \begin{pmatrix}0 & -1 & 0\\ 0 & 0 & 1\\ -1 & 0 & 0\end{pmatrix}.
  \end{matrix}
\end{equation}

\section{Finite-volume transport discretisation}

\label{sec:transport}

The transport scheme is an extension to
the method-of-lines scheme used by \citet{melvin19}. Solving the transport equation for an intermediate estimate $s^p$ of wave dynamics terms of a scalar field $s$ with a prescribed wind field $\mathbf{u}$ over a timestep either in advective form gives
\begin{equation}
    s^{n+1} = s^{p} - \Delta t\mathcal{A}\left(s^p,\mathbf{u}\right),\label{eq:transport-advective}
\end{equation}
or in flux form gives
\begin{equation}
    s^{n+1} = s^{p} - \Delta t \nabla\cdot\boldsymbol{\mathcal{F}}\left(s^p,\mathbf{u}\right).\label{eq:transport-flux}
\end{equation}
The transport scheme is chosen such that it maintains a number of desirable properties:
\begin{enumerate}
    \item The flexibility to alter the order of accuracy independently of the accuracy chosen for the finite-element, wave dynamics part of the model;
    \item At least second order temporal accuracy;
    \item Small dispersive errors;
    \item Scale selective damping;
    \item Stability for large CFL $\left(\mathbf{u}\Delta t/\Delta x>1\right)$ flows in any coordinate direction;
    \item Flux-form variables should return increments that are linear in the divergence;
    \item Computationally efficient.
\end{enumerate}

Achieving the first property rules out using the native finite-element discretisation when using degree $l=0$ spaces, where the transport scheme would be at best first-order, and instead motivates the use of a finite-volume method where the spatial accuracy can be linked to the polynomial reconstruction. The desired temporal accuracy is achieved through using an explicit multi-stage Runge-Kutta integration scheme for the transport terms in the same manner as \citet{melvin19}. The third and fourth properties can be achieved through using an upwind, even-degree polynomial reconstruction of the field used in the advective $\mathcal{A}$ and flux $\boldsymbol{\mathcal{F}}$ terms of \eqref{eq:transport-advective} and \eqref{eq:transport-flux} respectively. To achieve the fifth property the explicit Runge-Kutta scheme can be substepped with the transport scheme such that the effective CFL number for each substep is within the stability envelope of the desired Runge-Kutta scheme. The sixth property is motivated by a desire to maintain a constant density in non-divergent flows. To achieve this it is sufficient that the increment of a flux-form variable is linear in the divergence of the transporting wind field, here this is achieved by using an advective-then-flux scheme, Section \ref{sec:adv-then-flux}. 

As a step towards obtaining the seventh property the transport scheme is temporally split between the horizontal and vertical directions using a 2nd order Strang splitting (\citet{Strang68}). Splitting can deliver significant computational performance benefits, by allowing different schemes or options to be used for the different directions. Horizontal transport requires data communication costs when using multi-processor computers, while vertical transport typically involves higher Courant numbers (due to the grid anisotropy). Separating the two allows these problems to be addressed separately.

The horizontal spatial reconstruction follows that described in \citet{kent22} and the vertical reconstruction is extended from that used in \citet{melvin19} as described in Section \ref{sec:spatial-reconstruction}. The scheme defined in this section is applied in computational space, except for the spatial reconstructions which are computed in physical space. The finite-volume transport scheme used here is designed for a mesh with a single scalar degree of freedom per cell entity (i.e. in the cell centre or centre of a face). To couple the finite-element wave dynamics and the finite-volume transport, scalars need to be mapped from the finite-element spaces to the finite-volume space. The mapping from an degree $l=0$ finite-element space to the finite-volume scheme is the identity operator while mapping from an degree $l>0$ space requires a projection into the finite-volume space and is not considered here.

\subsection{Transported State}\label{sec:transport_state}

As in \citet{kent22} the transport acts on an intermediate state $\left(\mathbf{u}^p, \rho^p, \theta^p\right)$ (termed predictors in \citet{kent22}, their Section 4.3) for the prognostic variables $\left(\mathbf{u}, \rho, \theta\right)$ that consists of an explicit half timestep estimate of the wave dynamics terms only. The use of this estimate (instead of just taking the start of timestep fields) allows the present method to mimic the stability properties of a semi-implicit semi-Lagrangian scheme while using an Eulerian average advecting velocity $\overline{\mathbf{u}}^{1/2}$ rather rather than the Lagrangian average advecting velocity used in a standard semi-implicit semi-Lagrangian scheme. A full analysis and discussion of reasons behind this choice is in preparation and will be given elsewhere. The fields to be transported are given by
\begin{eqnarray}
\mathbf{u}^p &\equiv& \left[\mathbf{u} + \left(1-\alpha\right)\Delta t\mathbf{S}\right]^n,\label{eq:velocity-predictor}\\
\rho^p &\equiv& \left[\rho - \left(1-\alpha\right)\Delta t\rho\nabla\cdot\mathbf{u}\right]^n,\label{eq:continuity-predictor}\\
\theta^p &\equiv& \theta^n.\label{eq:thermo-predictor}
\end{eqnarray}
As in \citet{kent22} the choice of these terms is motivated by capturing the explicit parts of the non-advective processes e.g. $\rho^p$ contains the $\rho\nabla\cdot\mathbf{u}$ component of $\nabla\cdot\left(\rho\mathbf{u}\right)$ but not the advective $\mathbf{u}\cdot\nabla\rho$ component.

\subsection{Temporal Splitting: Advective Form}\label{sec:transport_strang}

The transport of a field $s$ by a three-dimensional velocity field $\mathbf{u}$ is split into vertical and horizontal components using a second-order Strang splitting as follows.
Noting that following \citet{mcrae16} (their table 3, where the $\wb_2$ space used here is referred to as $NCF_r$) the velocity space $\wb_2\equiv\wb_{2h}\oplus\wb_{2v}$ can be written as the composition of a space $\wb_2^h$ of vectors in the horizontal direction and $\wb_2^v$ of vectors in the vertical direction, (see \citet{maynard20}, their figure 1, for a illustrative example). The wind is split into horizontal $\mathbf{u}_H\in\wb_2^h$ and vertical $\mathbf{u}_V\in\wb_2^v$ components $\left(\mathbf{u} \equiv \mathbf{u}_H \oplus \mathbf{u}_V\right)$. The transport equation  $\partial s/\partial t + \mathcal{A}\left(s,\mathbf{u}\right)=0$ is then discretised across the timestep as
\begin{eqnarray}
    s^{\left(1\right)} &=& s^p - \frac{\Delta t}{2}\mathcal{A}\left(s^p,\mathbf{u}_V\right),\label{eq:strang1}\\
    s^{\left(2\right)} &=& s^{\left(1\right)} - \Delta t\mathcal{A}\left(s^{\left(1\right)},\mathbf{u}_H\right),\label{eq:strang2}\\
    s^{n+1} &=& s^{\left(2\right)} - \frac{\Delta t}{2}\mathcal{A}\left(s^{\left(2\right)},\mathbf{u}_V\right).\label{eq:strang3}
\end{eqnarray}
Each split step is then discretised using a multistage Runge-Kutta scheme as in  \citet{melvin19, kent22} which is further substepped to ensure stablility for large CFL numbers.

\subsection{Temporal Splitting: Flux Form}\label{sec:adv-then-flux}

It is tempting to formulate the flux form by replacing the $\mathcal{A}$'s in \eqref{eq:strang1}-\eqref{eq:strang3} with $\nabla\cdot\boldsymbol{\mathcal{F}}$'s. However, in the presence of non-divergent flow this does not have the desirable property of preserving a constant (except for the trivial case where the components of the divergence in each direction are zero) as $s^{n+1}$ is not linear in $\nabla\cdot\mathbf{u}$. Simply replacing only the $\mathcal{A}$ in \eqref{eq:strang3} by $\nabla\cdot\boldsymbol{\mathcal{F}}$ does not conserve the total scalar being transported. However, the following split form is both conservative and preserves a constant (\citet{bendall23}).

Revisiting the Strang splitting, the final step \eqref{eq:strang3} is modified as follows 
\begin{equation}
    s^{n+1} = s^p - \frac{\Delta t}{2}\nabla\cdot\boldsymbol{\mathcal{F}}\left(s^p,\mathbf{u}_V\right) - \Delta t\nabla\cdot\boldsymbol{\mathcal{F}}\left(s^{\left(1\right)},\mathbf{u}_H\right) - \frac{\Delta t}{2}\nabla\cdot\boldsymbol{\mathcal{F}}\left(s^{\left(2\right)},\mathbf{u}_V\right).\label{eq:advective-then-flux3}
\end{equation}
When $s^p$ is constant the final step is seen to reduce to $s^{n+1} = s^p\left(1 - \Delta t\nabla\cdot\mathbf{u}\right)$ which preserves $s^p$ when $\nabla\cdot\mathbf{u}=0$.

At any stage of the split step, where a flux-form equation is solved, such as \eqref{eq:advective-then-flux3}, the advective-then-flux form of a Runge-Kutta scheme can be obtained for an $m-$stage scheme as
\begin{eqnarray}
    s^{\left(i\right)} &=& s^p - \Delta t\sum_{j=1}^{i-1}a_{i,j}\mathcal{A}\left(s^{\left(j\right)},\mathbf{u}\right),\quad i=1,\ldots,m,\label{eq:RK-advective-then-flux1}\\
    s^{n+1} &=& s^p - \Delta t\sum_{k=1}^{m}b_k\nabla\cdot\boldsymbol{\mathcal{F}}\left(s^{\left(k\right)},\mathbf{u}\right),
\end{eqnarray}
where the coefficients $a_{i,j}$ and $b_k$ are given by the Butcher tableau (\citet{butcher87}) associated with the chosen scheme and again when $s^p$ is constant the final step reduces to $s^{n+1} = s^p\left(1 - \Delta t\nabla\cdot\mathbf{u}\right)$.

\subsection{Spatial reconstruction of a scalar field}\label{sec:spatial-reconstruction}

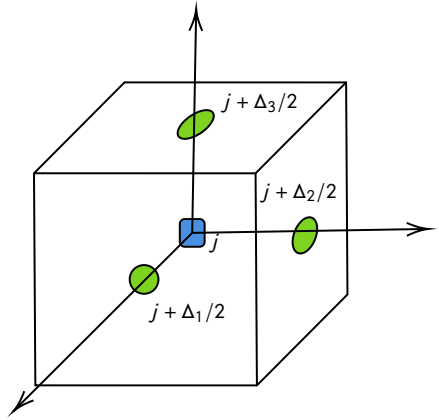
\begin{figure}[htbp]
\centering
\tikzset{every picture/.style={line width=0.75pt}} 

\begin{tikzpicture}[x=0.75pt,y=0.75pt,yscale=-1,xscale=1]

\draw   (264.81,109.89) -- (310.41,63.99) -- (422.01,64.01) -- (422.69,171.11) -- (377.09,217.01) -- (265.49,216.99) -- cycle ; \draw   (422.01,64.01) -- (376.41,109.91) -- (264.81,109.89) ; \draw   (376.41,109.91) -- (377.09,217.01) ;
\draw  [fill={rgb, 255:red, 74; green, 144; blue, 226 }  ,fill opacity=1 ] (338.25,135.5) .. controls (338.25,134.12) and (339.37,133) .. (340.75,133) -- (348.25,133) .. controls (349.63,133) and (350.75,134.12) .. (350.75,135.5) -- (350.75,144.5) .. controls (350.75,145.88) and (349.63,147) .. (348.25,147) -- (340.75,147) .. controls (339.37,147) and (338.25,145.88) .. (338.25,144.5) -- cycle ;
\draw  [fill={rgb, 255:red, 126; green, 211; blue, 33 }  ,fill opacity=1 ] (313,163.5) .. controls (313,159.5) and (316.25,156.25) .. (320.25,156.25) .. controls (324.25,156.25) and (327.5,159.5) .. (327.5,163.5) .. controls (327.5,167.5) and (324.25,170.75) .. (320.25,170.75) .. controls (316.25,170.75) and (313,167.5) .. (313,163.5) -- cycle ;
\draw  [fill={rgb, 255:red, 126; green, 211; blue, 33 }  ,fill opacity=1 ] (407.43,137.46) .. controls (407.6,142.05) and (404.97,147.47) .. (401.56,149.57) .. controls (398.15,151.66) and (395.24,149.64) .. (395.07,145.04) .. controls (394.9,140.45) and (397.53,135.03) .. (400.94,132.93) .. controls (404.35,130.84) and (407.26,132.86) .. (407.43,137.46) -- cycle ;
\draw  [fill={rgb, 255:red, 126; green, 211; blue, 33 }  ,fill opacity=1 ] (339,85.25) .. controls (342.04,81.25) and (347.75,78) .. (351.75,78) .. controls (355.75,78) and (356.54,81.25) .. (353.5,85.25) .. controls (350.46,89.25) and (344.75,92.5) .. (340.75,92.5) .. controls (336.75,92.5) and (335.96,89.25) .. (339,85.25) -- cycle ;
\draw    (344.5,140) -- (461.5,138.03) ;
\draw [shift={(463.5,138)}, rotate = 179.04] [color={rgb, 255:red, 0; green, 0; blue, 0 }  ][line width=0.75]    (10.93,-3.29) .. controls (6.95,-1.4) and (3.31,-0.3) .. (0,0) .. controls (3.31,0.3) and (6.95,1.4) .. (10.93,3.29)   ;
\draw    (344.5,140) -- (346.46,29) ;
\draw [shift={(346.5,27)}, rotate = 91.01] [color={rgb, 255:red, 0; green, 0; blue, 0 }  ][line width=0.75]    (10.93,-3.29) .. controls (6.95,-1.4) and (3.31,-0.3) .. (0,0) .. controls (3.31,0.3) and (6.95,1.4) .. (10.93,3.29)   ;
\draw    (344.5,140) -- (255.91,228.59) ;
\draw [shift={(254.5,230)}, rotate = 315] [color={rgb, 255:red, 0; green, 0; blue, 0 }  ][line width=0.75]    (10.93,-3.29) .. controls (6.95,-1.4) and (3.31,-0.3) .. (0,0) .. controls (3.31,0.3) and (6.95,1.4) .. (10.93,3.29)   ;

\draw (352.75,138.9) node [anchor=north west][inner sep=0.75pt]    {$j$};
\draw (378.41,113.31) node [anchor=north west][inner sep=0.75pt]    {$j+\Delta _{2}/2$};
\draw (322.25,174.15) node [anchor=north west][inner sep=0.75pt]    {$j+\Delta _{1}/2$};
\draw (358.75,67.9) node [anchor=north west][inner sep=0.75pt]    {$j+\Delta _{3}/2$};

\end{tikzpicture}
\begin{center}

\end{center}
\caption{Cell $j$ with points staggered half a grid point located on cell faces $\Delta_i/2,$ $i=1,2,3$. For a field stored at $j$ the spatial reconstruction is computed at $j\pm\Delta_i/2$.}
\label{fig:cube}
\end{figure}

The transport scheme computes a high order reconstruction $\check{s}$ of a given scalar field $s$ in physical space $\boldsymbol{\chi}$. For a scalar value in cell $j$ the reconstructed field is computed at points staggered half a grid length $j\pm \Delta_i/2, \quad i=1,2,3$ from the original field in all three directions, where $\Delta_i$ is defined to be the grid spacing in the $\widehat{\chi}_i$ direction. For example, for a field in $\wb_3$ such as the density $\rho_j$, which is located at cell centres, then the reconstructed field $\check{\rho}_{j\pm\Delta_i/2}$ is computed at the centre of each cell face $j\pm\Delta_i/2$, see Figure \ref{fig:cube}. The reconstruction is computed by fitting a polynomial through a number of cells and evaluating this polynomial at the staggered points. The reconstruction is given an upwind bias, determined by the wind direction ($\mathbf{u}_j$ for $\mathcal{A}$ and $\mathbf{u}_{j\pm\delta_i}$ for $\boldsymbol{\mathcal{F}}$), by using even order polynomials for the reconstruction. 

%

The horizontal spatial reconstruction is the same as \citet{kent22} and is based on that used in \citet{thuburncottdubos13} and the interested reader is referred to \citet{baldauf08} and \citet{skammenc10} for other results on these types of schemes. In brief, a two-dimensional polynomial in local Cartesian coordinates is fitted in a least squares sense to a region of cells around the reconstruction point. This polynomial is then evaluated at the reconstruction point to give the reconstructed field. This method results in a scheme that is accurate across discontinuities in the mesh (such as at panel boundaries) and reduces grid imprinting from the transport scheme.


The vertical reconstruction follows the same method as the horizontal scheme, except now a one-dimensional polynomial is used and the local coordinate system $z$ can be taken to be aligned to the radial distance from the surface of the sphere and the origin is the at the reconstruction point. 

Near the top and bottom boundaries when there are not enough points to construct an upwind degree $n$ polynomial then the stencil is shifted so that the same polynomial is used as the first point where there are enough points for an upwind polynomial. The result is that the polynomial is no longer upwinded but the desired degree is maintained.

\subsection{Flux computation}\label{sec:flux_comp}
The flux $\boldsymbol{\mathcal{F}}$ is computed as in \citet{kent22} by a pointwise multiplication of the computational wind field $\widehat{\mathbf{u}}$ sampled at $j\pm\delta_i$ by the reconstructed scalar $\check{s}$
\begin{equation}
    \boldsymbol{\mathcal{F}} \equiv \widehat{\mathbf{u}}\left(\widehat{\boldsymbol{\chi}}_{j\pm\delta_i}\right)\check{s}. \label{eq:flux}
\end{equation}
The finite-volume divergence operator for a flux is then given by
\begin{equation}
    \nabla\cdot\boldsymbol{\mathcal{F}} = \frac{1}{\detj}\sum_{j=1,\ldots,6} \boldsymbol{\mathcal{F}}\left(\widehat{\boldsymbol{\chi}}_j\right)\cdot \widehat{\mathbf{n}}_j,
\end{equation}
where $j$ is the index of each face of the cell.

\subsection{Advective increment computation}\label{sec:adv_inc_comp}

To compute the advective tendency $\mathcal{A} \equiv \mathbf{u}\cdot\nabla s$ of a field $s$, given the reconstructed field $\check{s}$ located at points staggered half a grid length from $s_j$ the advective update is 
\begin{equation}
\widehat{\mathcal{A}}_j \equiv \frac{1}{\detj}\widehat{\mathbf{u}}\left(\widehat{\boldsymbol{\chi}}_j\right)\cdot \widehat{\boldsymbol{\delta}} \check{s}_j,\label{eq:MoL-advective-update}
\end{equation}
where $\widehat{\mathbf{u}}\left(\widehat{\boldsymbol{\chi}}_j\right)$ is the computational velocity field sampled at $\boldsymbol{\chi}_j$ and $\widehat{\boldsymbol{\delta}}$ is the discrete gradient operator in the computational space. For example if $\check{s}$ are located on cell faces  $\left[\check{s}_{j+\Delta_1/2},\check{s}_{j-\Delta_1/2},\check{s}_{j+\Delta_2/2},\check{s}_{j-\Delta_2/2},\check{s}_{j+\Delta_3/2},\check{s}_{j-\Delta_3/2}\right]$ on the (East, West, North, South, Up, Down) sides of cell $j$ respectively then 
\begin{equation}
\delta\check s_j \equiv \left[\left(\check{s}_{j+\Delta_1/2} - \check{s}_{j-\Delta_1/2}\right), 
  \left(\check{s}_{j+\Delta_2/2}-\check{s}_{j+\Delta_2/2}\right), \left(\check{s}_{j+\Delta_3/2}-\check{s}_{j-\Delta_3/2}\right)\right].
\end{equation}

Optionally monotonicity can be enforced on the advective update $\widehat{\mathcal{A}}$ through a simple clipping scheme. On the final stage of the Runge-Kutta scheme the update $\widehat{\mathcal{A}}$ is modified to ensure that $s^{min}\leq s^{n+1} \leq s^{max}$ where $s^{min}$ and $s^{max}$ are the minimum and maximum values of $s^p$ used in the stencil to compute $\widehat{\mathcal{A}}$ respectively. 

\subsection{Advection of vector fields}

To compute the advective increment $\widehat{\boldsymbol{\mathcal{A}}}\left(\widehat{\boldsymbol{\nu}},\uadvc\right)$ of a vector field $\widehat{\boldsymbol{\nu}}$ the advected field is first transformed into physical space using the $\wb_2$ transform 
\begin{equation}
\dot{\boldsymbol{\chi}}\equiv \frac{\partial\boldsymbol{\chi}}{\partial t} = \frac{\mathbf{J}\widehat{\boldsymbol{\nu}}}{\detj},
\end{equation}
where each component of $\dot{\boldsymbol{\chi}}$ is placed in the $\wb_3$ space $\dot{\chi}_i\in\wb_3\,i=1,2,3$. Since $\boldsymbol{\chi}$ is a Cartesian coordinate system this means that the Cartesian components of the velocity vector are transported. Each component of the vector is transported as a scalar in the advective form by solving
\begin{equation}
    \dot{\chi}_i^{n+1} = \dot{\chi}^p_i - \mathcal{A}\left(\dot{\chi}_i,\uadvc\right),\label{eq:chi-advection}
\end{equation}
using the method described in the previous sections. Once all three components have been advected the advective increment $\widehat{\boldsymbol{\mathcal{A}}}$ is
\begin{equation}
    \widehat{\boldsymbol{\mathcal{A}}} \equiv \left[\widehat{\mathcal{A}}\left(\dot{\chi}_1,\uadvc\right),\widehat{\mathcal{A}}\left(\dot{\chi}_2,\uadvc\right),\widehat{\mathcal{A}}\left(\dot{\chi}_3,\uadvc\right)\right],\label{eq:velocity-advective-increment}
\end{equation}
and is mapped back to the $\wb_2$ space through \eqref{eq:rhs-u-adv}.

\section{Timestepping}\label{sec:timestepping}

The timestepping algorithm closely follows that of \citet{melvin19} but, 
inspired by the algorithm of \citet{wood14}, the Newton loop is split into $(n_o)$ outer loops and $(n_i)$ inner loops, such that the advection scheme is called $n_o$ times in the outer loop, whilst updates to the residuals and the linear solver are called $n_o\times n_i$ times in the inner loop. 

\subsection{Linear System}

\label{sec:linear-system}

The solution procedure follows \citet{melvin19} except that the linear system
 is modified to include the Coriolis terms. To recap, the linear system
\begin{equation}
{\cal L}\left(\mathbf{x}^{*}\right)\mathbf{x}^{\prime}=-{\cal R}\left(\mathbf{x}^{\left(k\right)}\right),\quad \mathbf{x}\equiv\left[\mathbf{u}, \rho, \theta, \Pi\right]^T,\label{eq:quasi-Newton}\end{equation}
is solved, with the linear operator ${\cal L}$ inspired by the linearisation of the set of residuals $\mathcal{R}$ (see section \ref{sec:discrete-linear-system} below) about some reference state $\mathbf{x}^{*}\equiv\left[\boldsymbol{0}, \rho^*, \theta^*, \Pi^*\right]^T$ with relaxation factors $\tau_{u,\rho,\theta}$ to obtain
${\cal L}\left(\mathbf{x}^{*}\right)$.
In spatially continuous form $\mathcal{L}$ from \citet{melvin19} (their (42)) is augmented by the
Coriolis terms to give
\begin{equation}
{\cal L}\left(\mathbf{x}^{*}\right)\mathbf{x}^{\prime}=\begin{cases}
\left(1+ 2\tau_u\Delta t \boldsymbol{\Omega}\times\right)\mathbf{u}^{\prime}-\mu\left(\frac{\mathbf{n}_{b}\cdot\mathbf{u}^{\prime}}{\mathbf{n}_{b}\cdot\mathbf{z}_{b}}\right)\mathbf{z}_{b}\\
\,\,\,\,\,+\tau_{u}\Delta tc_{p}\left(\theta^{\prime}\mathbf{z}_{b}\left[\mathbf{n}_b\cdot\nabla\Pi^{*}\right]+\theta^{*}\nabla\Pi^{\prime}\right),\\
\rho^{\prime}+\tau_{\rho}\Delta t\nabla\cdot\left(\rho^{*}\mathbf{u}^{\prime}\right),\\
\theta^{\prime}+\tau_{\theta}\Delta t\mathbf{u}^{\prime}\cdot\mathbf{z}_{b}\left(\mathbf{n}_b\cdot\nabla\theta^{*}\right),\\
\frac{1-\kappa}{\kappa}\frac{\Pi^{\prime}}{\Pi^{*}}-\frac{\rho^{\prime}}{\rho^{*}}-\frac{\theta^{\prime}}{\theta^{*}}.\end{cases}\label{eq:L}\end{equation}
Note, following \citet{wood14} but in contrast to \citet{melvin19}, only the vertical component of the implicit buoyancy terms $\left(\theta'\nabla\Pi^*,\,\mathbf{u}'\cdot\nabla\theta^*\right)$ are retained in the linear system which improves the efficiency of the solver without detracting from the convergence properties.

\subsection{Inner Loop Convergence}
\label{sec:inner-loop}

The implicit terms in the continuity and thermodynamic equations contain corrections to the transport terms, $\nabla\cdot\left(\rho^*\mathbf{u}'\right)$ and $\mathbf{u}'\cdot\nabla\theta^*$ respectively, where the velocity increment $\mathbf{u}'\equiv\mathbf{u}^{\left(k+1\right)}-\mathbf{u}^{\left(k\right)}$ is defined relative to the latest (inner loop) estimate for the velocity field $\mathbf{u}^{\left(k\right)}$. Since the transport terms are only updated in the outer loop the advecting wind is $\overline{\mathbf{u}}^{1/2} = 1/2\left(\mathbf{u}^{\left(o\right)}+\mathbf{u}^{n}\right)$ where $\mathbf{u}^{\left(o\right)}$ is the latest estimate of $\mathbf{u}^{n+1}$ available in the outer loop. If the inner loop iterate $n_i>1$ then the two estimates for $\mathbf{u}^{n+1}$ do not agree $\left(\mathbf{u}^{\left(o\right)}\neq \mathbf{u}^{\left(k\right)}\right)$ and this can lead to an inconsistency in the discretisation. Taking the thermodynamic equation $\theta' + \tau_\theta\Delta t\mathbf{u}'\cdot\nabla\theta^* = R_\theta$ with $\tau_\theta=1/2$ and $\theta^* = \theta^n$ this becomes
\begin{eqnarray}
    \theta^{\left(k+1\right)} &-& \theta^{\left(k\right)} + \frac{\Delta t}{2}\left(\mathbf{u}^{\left(k+1\right)}-\mathbf{u}^{\left(k\right)}\right)\cdot\nabla\theta^n \nonumber\\ &=& R_\theta \equiv -\left[\theta^{\left(k\right)} - \theta^n + \frac{\Delta t}{2}\left(\mathbf{u}^{\left(o\right)}+\mathbf{u}^{n}\right)\cdot\nabla\theta^n\right],\label{eq:inner-loop-problem1}
\end{eqnarray}
which can be rearranged to give
\begin{eqnarray}
    \theta^{\left(k+1\right)} - \theta^n &+& \frac{\Delta t}{2}\left(\mathbf{u}^{\left(k+1\right)}+\mathbf{u}^{n}\right)\cdot\nabla\theta^n \nonumber\\&=& \frac{\Delta t}{2}\left(\mathbf{u}^{\left(k\right)}-\mathbf{u}^{\left(o\right)}\right)\cdot\nabla\theta^n, \label{eq:inner-loop-problem2}
\end{eqnarray}
where the left hand side is the desired temporal discretisation of the equation (centred implicit). However the term on the right hand side only vanishes if $\mathbf{u}^{\left(k\right)}=\mathbf{u}^{\left(o\right)}$. The solution to this inconsistency is that in the inner loop when $n_i>1$ the residuals of the thermodynamic $\left(R_\theta\right)$ and continuity $\left(R_\rho\right)$ equations are set to zero. Returning to the above example \eqref{eq:inner-loop-problem1} is replaced by
\begin{equation}
  \theta^{\left(k+1\right)} - \theta^{\left(k\right)} + \frac{\Delta t}{2}\left(\mathbf{u}^{\left(k+1\right)}-\mathbf{u}^{\left(k\right)}\right)\cdot\nabla\theta^n = 0, \label{eq:inner-loop-solution1}
\end{equation}
and for $k$ such that $\mathbf{u}^{\left(k-1\right)}=\mathbf{u}^{\left(o\right)}$ then replacing $\theta^{\left(k\right)}$ with \eqref{eq:inner-loop-problem1} with all indices reduced by $1$ yields the desired form
\begin{equation}
    \theta^{\left(k+1\right)} - \theta^n + \frac{\Delta t}{2}\left(\mathbf{u}^{\left(k+1\right)}+\mathbf{u}^{n}\right)\cdot\nabla\theta^n = 0. \label{eq:inner-loop-solution2}
\end{equation}
The same result then applies for larger $k$ upon repeated use of \eqref{eq:inner-loop-solution2}.
This inner loop correction to the residuals is not applied to the momentum equation or equation of state since in this case the linear corrections in \eqref{eq:L} correspond to terms that are updated in the inner loop, i.e. there are no linear transport terms, and so there is no inconsistency.

\subsection{Discrete Linear System}
\label{sec:discrete-linear-system}

Applying
the mixed finite-element discretisation of Section \ref{sec:Spatial-discretization} to \eqref{eq:L} results in:\begin{eqnarray}
M_{2}^{\mu,C}\widetilde{u}^{\prime}-P_{2v\theta}^{\Pi^{*}}\widetilde{\theta}^{\prime}-G^{\theta^{*}}\widetilde{\Pi}^{\prime} & = & -{\cal R}_{\mathbf{u}},\label{eq:si_u_discrete}\\
M_{3}\widetilde{\rho}^{\prime}+D\left(\widehat{\rho}^{*}\widetilde{u}^{\prime}\right) & = & -{\cal R}_{\rho},\label{eq:si_rho_discrete}\\
M_{\theta}\widetilde{\theta}^{\prime}+P_{\theta2v}^{\theta^{*}}\widetilde{u}^{\prime} & = & -{\cal R}_{\theta},\label{eq:si_theta_discrete}\\
E^{\Pi^{*}}\widetilde{\Pi}^{\prime}-E^{\rho^{*}}\widetilde{\rho}^{\prime}-E^{\theta^{*}}\widetilde{\theta}^{\prime} & = & -{\cal R}_{\Pi},\label{eq:si_exner_discrete}\end{eqnarray}
with
\begin{eqnarray}
{\cal R}_{\mathbf{u}} &=& \Delta t\left(M_2\delta_t\widetilde{u} + M_\mu\overline{\widetilde{u}}^1 - \overline{R_u}^\alpha - R_u^A\right),\label{eq:RHS_u}\\
{\cal R}_{\rho} &=& \Delta t\left(M_3\delta_t\widetilde{\rho} - R_\rho^F\right),\label{eq:RHS_theta}\\
{\cal R}_{\theta} &=& \Delta t\left(M_\theta\delta_t\widetilde{\theta} - R_\theta^A\right),\label{eq:RHS_rho}\\
{\cal R}_{\Pi} &\equiv& \left(1-\frac{p_0\Pi^{\frac{1-\kappa}{\kappa}}}{R\rho\theta}\right),\label{eq:RHS_exner}
\end{eqnarray}
 where $P_{2v\theta}^{\Pi^*}$ and $P_{\theta 2v}^{\theta^*}$ are the vertically restricted versions of the $P_{2\theta}^{\Pi^*}$ and $P_{\theta2}^{\theta^*}$ operators given in \citet{melvin19} (their (84) and (87)). $M_{2}^{\mu,C}$ is the operator formed by combining
the $\mathbb{W}_{2}$ mass matrix with the operators arising from the
Rayleigh damping and Coriolis terms:\begin{equation}
M_{2}^{\mu,C}\equiv M_{2}+\Delta tM_{\mu}+\tau_u\Delta tM_{C},\label{eq:M2mu}\end{equation}
with \begin{equation}
\left(M_{C}\right)_{ij}\equiv\left\langle \mathbf{J}\widehat{\mathbf{v}}_{i},2\boldsymbol{\Omega}\times\frac{\mathbf{J}\widehat{\mathbf{v}}_{j}}{\det\left(\mathbf{J}\right)}\right\rangle.\label{eq:MC}\end{equation}
As, in contrast to \citet{melvin19}, the equation of state is now sampled the operators in \eqref{eq:si_exner_discrete} are given by
\begin{eqnarray}
  E^{\Pi^*} &=& \frac{1-\kappa}{\kappa}\left[\frac{p_0}{R}\frac{\left(\widehat{\Pi^*}\right)^{\frac{1-\kappa}{\kappa}}}{\widehat{\rho^*}\widehat{\theta^*}}\right]\frac{\widehat{\sigma}}{\widehat{\Pi^*}},\label{eq:E_pi}\\
  E^{\rho^*} &=& \frac{\widehat{\sigma}}{\widehat{\rho^*}},\label{eq:E_rho}\\
  E^{\theta^*} &=& \frac{\widehat{w}}{\widehat{\theta^*}},\label{eq:E_theta}
\end{eqnarray}
where each entry $E_{i,j}$ of an operator $E$ is obtained by evaluating all variables at nodal point $\widehat{\boldsymbol{\chi}_i}$ with basis function $\widehat{\sigma}_j$ or  $\widehat{w}_j$.
All other operators have the same form as given by \citet{melvin19}.
At convergence of the iterative procedure primed quantities vanish and ${\cal R}\left(\mathbf{x}^{\left(k\right)}\right)=0$ is solved.

\subsection{Iterative Solver}
\label{sec:iterative-solver}

The system of equations \eqref{eq:si_u_discrete}-\eqref{eq:si_exner_discrete}
is solved using the method presented in \citet{maynard20}. This consists of  an iterative Krylov method that is preconditioned by an approximate Schur complement of the equations for the pressure increment which is itself solved using a single v-cycle of a geometric multigrid method. The approximate Schur complement is achieved by a diagonal mass lumping of $M_\theta$ and $M_2^{\mu}$ (where the Coriolis terms have been dropped from the lumped approximation, equivalent to $\tau_u=0$ in \eqref{eq:M2mu}).

\section{Computational examples}

\label{sec:Results}

In order to assess the accuracy of the model it is run on a set of standard numerical tests for atmospheric dynamics drawn from the literature. These are used to ensure the model generates the correct response to forcing at different scales as well as maintaining the large scale balances important in the governing equations. Of particular concern with meshes, such as the cubed-sphere, that have discontinuities in their coordinate lines, is what impact those discontinuities have on the numerical solutions, referred to as grid imprinting.

The tests results presented in this section are:
\begin{description}
  \item[\ref{sec:resting atmosphere}] Resting atmosphere over orography (Test 2.0.0 from the 2012 DCMIP project, \citet{dcmip2012})
  \item[\ref{sec:gaussian}] Flow over a Gaussian mountain (\citet{allen16})
  \item[\ref{sec:baroclinic}] Deep atmosphere baroclinic wave (\citet{ullrich14})
  \item[\ref{sec:held-suarez}] Held-Suarez climate test (\citet{HeldSuarez94})
\end{description}
Key test parameters for each of these examples are summarised in Table \ref{tab:Test_params-1} where the average grid spacing has been taken to be the square root of the average cell area. For a $Cn$ mesh this is given by
\begin{equation}
    dC = \sqrt{\frac{4\pi a^2}{6n^2}}. \label{eq:average-gridspacing}
\end{equation}

While in principle the finite-element
methodology affords flexibility in the polynomial degree, as in \citet{melvin19} the focus will again be on results in the lowest-degree $l=0$ case. Additionally a number of simplifications and specifications to the formulation given in previous sections are made:
\begin{itemize}
\item The coordinate space is $\wb_\chi = Q_1^{DG}$; 
\item The angular resolution of the cubed-sphere mesh is kept constant: $\Delta\xi_c = \Delta \eta_c = \frac{\pi}{2n},\quad \forall c$;
\item The semi-implicit scheme is centred in time: $\alpha=1/2$;
\item The relaxation parameters are $\tau_{u}=1/2$ and $\tau_{\rho,\theta}=1$ which is empirically found to improve convergence, consistent with \citet{wood14};
\item 2 outer (advection) and 2 inner (nonlinear) iterations are used: $n_o = n_i = 2$;
\item As in \citet{melvin19} the reference profiles $\mathbf{x}^*$ are taken to be the start of timestep fields $\mathbf{x}^*\equiv \mathbf{x}^n$ with no adjustment applied to these fields;
\item A quadratic reconstruction of scalar fields is used in the advection scheme to compute  fluxes $\widehat{\boldsymbol{\mathcal{F}}}$ and advective updates $\widehat{\mathcal{A}}$;
\item The damping layer is not required for any tests considered here and so $\mu\equiv 0$;
\item The vertical mesh consists of $n$ levels and the height of each level $k$ is
\begin{equation}
  z_k = z_T\epsilon_k + z_B\left(1-\epsilon_k\right)
\end{equation}
where $z=r-a$ is the height above the surface of the sphere and $z_B$ is height of the domain surface above $a$ ($z_B=0$ or as given by any orographic profile). The non-dimensional parameter $\epsilon$ is given by
\begin{equation}
   \epsilon_k = \left(k/n\right),\label{eq:uniform-mesh}
\end{equation}
for a uniform vertical mesh or 
\begin{equation}
    \epsilon_k = \frac{\sqrt{\gamma\left(k/n\right)^2+1}-1}{\sqrt{\gamma + 1}-1},\label{eq:quadratic-mesh}
\end{equation}
with $\gamma = 15$ for a quadratic stretching.
\item All initial conditions are computed by sampling the field at degree of freedom locations and where required the density or initial pressure are obtained from the equation of state. Additionally, no discrete balance is applied to the initial conditions.
\end{itemize}
For the spherical domain used in these examples the coordinate $\boldsymbol{\chi}$ is replaced by the standard geocentric Cartesian coordinates $\mathbf{X}$ so that $\left(\chi_1,\,\chi_2,\,\chi_3 \right) \equiv \left(X,\,Y,\,Z \right)$. Additionally the results are linearly interpolated into a regular lat-long grid for presentation. In this section $w$ is used to denote the vertical component of the velocity $\mathbf{u}$ in the radial direction (i.e.~$w=Dr/Dt\neq D\chi_3/Dt$) and $u$ is the zonal component of the velocity $\mathbf{u}$ (i.e ~$u=r\cos\phi D\lambda/Dt\neq D\chi_1/Dt$). Initial conditions are given in spherical coordinates with latitude $\lambda\in\left(-\pi/2,\pi/2\right)$ and longitude $\phi\in\left(-\pi,\pi\right)$.

\begin{table*}
\begin{centering}
\begin{tabular}{|c|c|c|c|c|c|c|}
\hline
Test & Section & Resolution $CnLm$ & Approx grid &  Model depth & $\Delta t$ \tabularnewline
 & &  ($n$ cells per edge, &  spacing $dC$ (km) & $z_T$ (km) & (s) \tabularnewline
  & &  $m$ layers) &   &  & \tabularnewline
\hline
Resting atmosphere & \ref{sec:resting atmosphere} & C96L30 & 96.0 & 12 & 600 \tabularnewline
Gaussian mountain & \ref{sec:gaussian} & C96L40 & 96.0 &  32 & 900 \tabularnewline
Baroclinic wave & \ref{sec:baroclinic} & C96L30 & 96.0 &  30 & 900 \tabularnewline
Baroclinic wave: grid imprinting & \ref{sec:hiwpp} & C448L30 & 20.6 & 30 & 225 \tabularnewline
 &  & C896L30 & 10.3 & 30 & 225 \tabularnewline
Held-Suarez & \ref{sec:held-suarez} & C48L30 & 192.1 & 30 & 1800 \tabularnewline
\hline
\end{tabular}
\par\end{centering}
\caption{\label{tab:Test_params-1}Model parameters for each test.}\end{table*}

\subsection{Resting atmosphere over orography}
\label{sec:resting atmosphere}

Orography is represented in the model formulation through the Piola transforms and in particular through the Jacobian $\mathbf{J}$. In the presence of orography the mapping $\partial\boldsymbol{\xi}/\partial\widehat{\boldsymbol{\chi}}$ \eqref{eq:dxidchi} introduces a coupling between the horizontal components of the velocity vector and the vertical component of the momentum equation (through the $\tan\delta_\xi$ and $\tan\delta_\eta$ terms in \eqref{eq:dxidchi}). This is in contrast to most models (such as \citet{wood14}) where the presence of orography and terrain following coordinates introduces a coupling of the vertical components of the pressure gradient term into the horizontal components of the momentum equation  (for example (22) \& (23) of \citet{wood14}). It is therefore interesting to see how this different formulation can represent a balanced state over orography. To complement the orographic tests already presented in \citet{melvin19}  test 2.0.0 from the DCMIP2012
project \citet{dcmip2012} is used to simulate a resting atmosphere over large scale orography.

\begin{figure}[htbp]
\centering
\subfloat[SI-FE/FV $u$]{\includegraphics[scale=0.3]{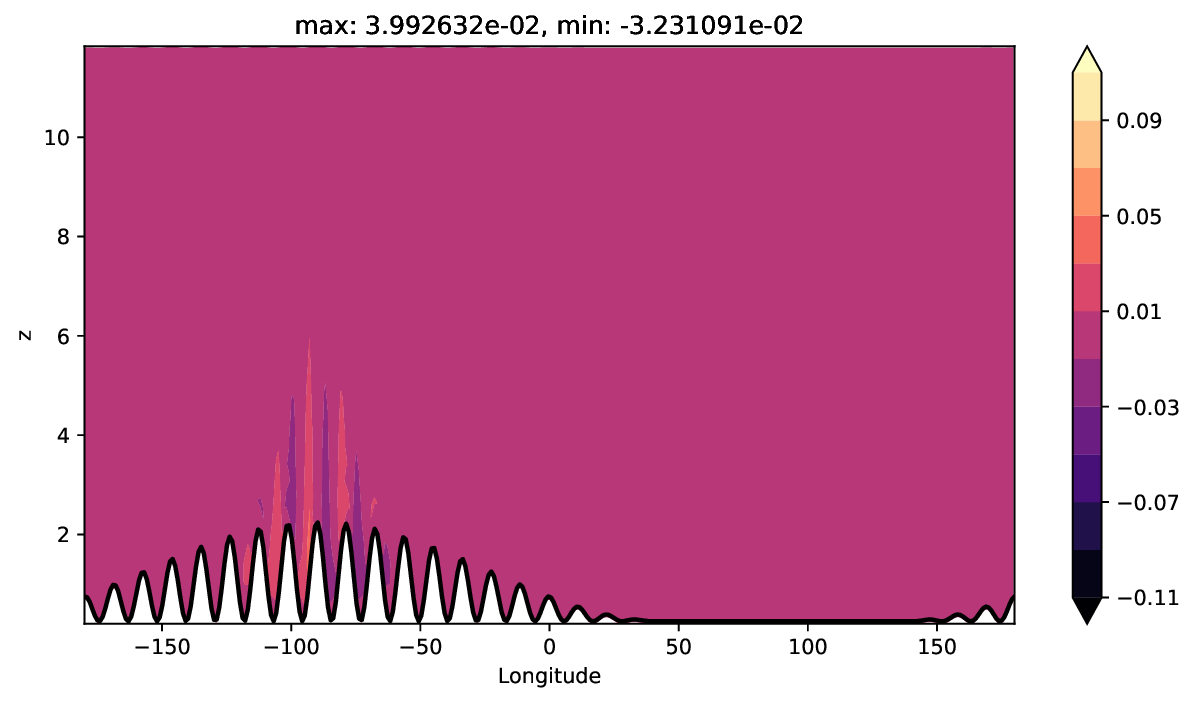}}%
\subfloat[SI-FE/FV $w$]{\includegraphics[scale=0.3]{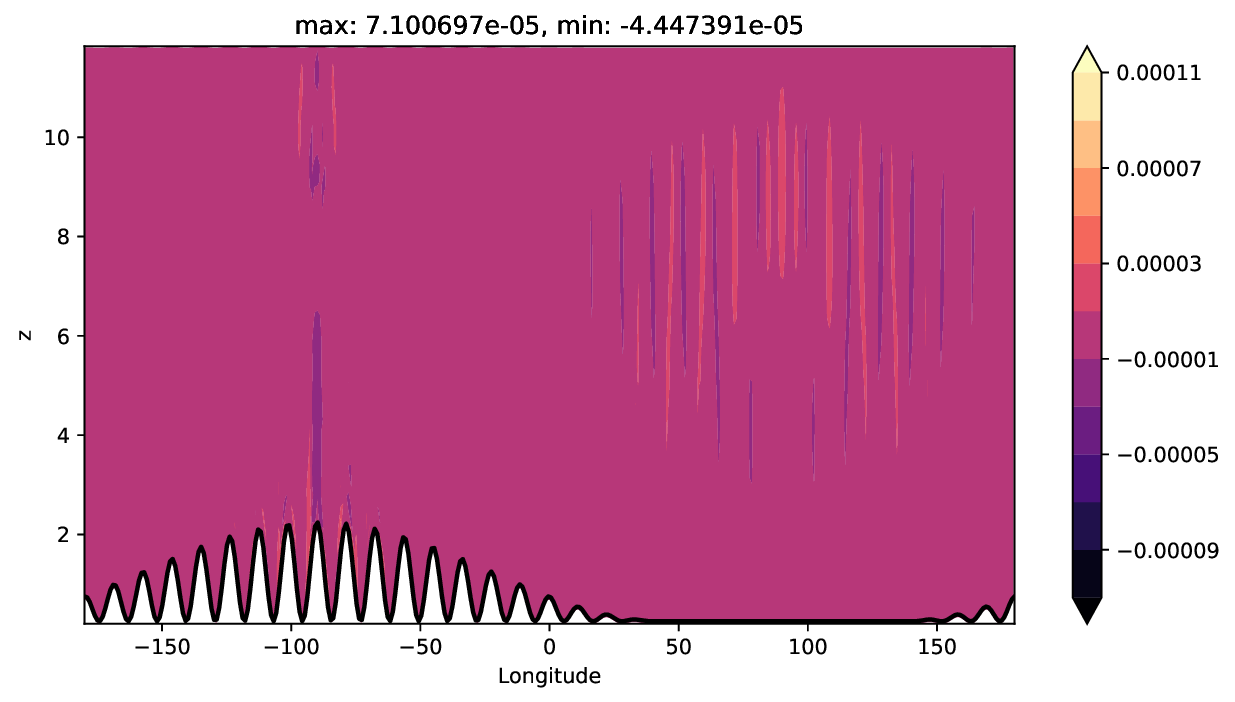}}\\
\subfloat[SISL $u$]{\includegraphics[scale=0.3]{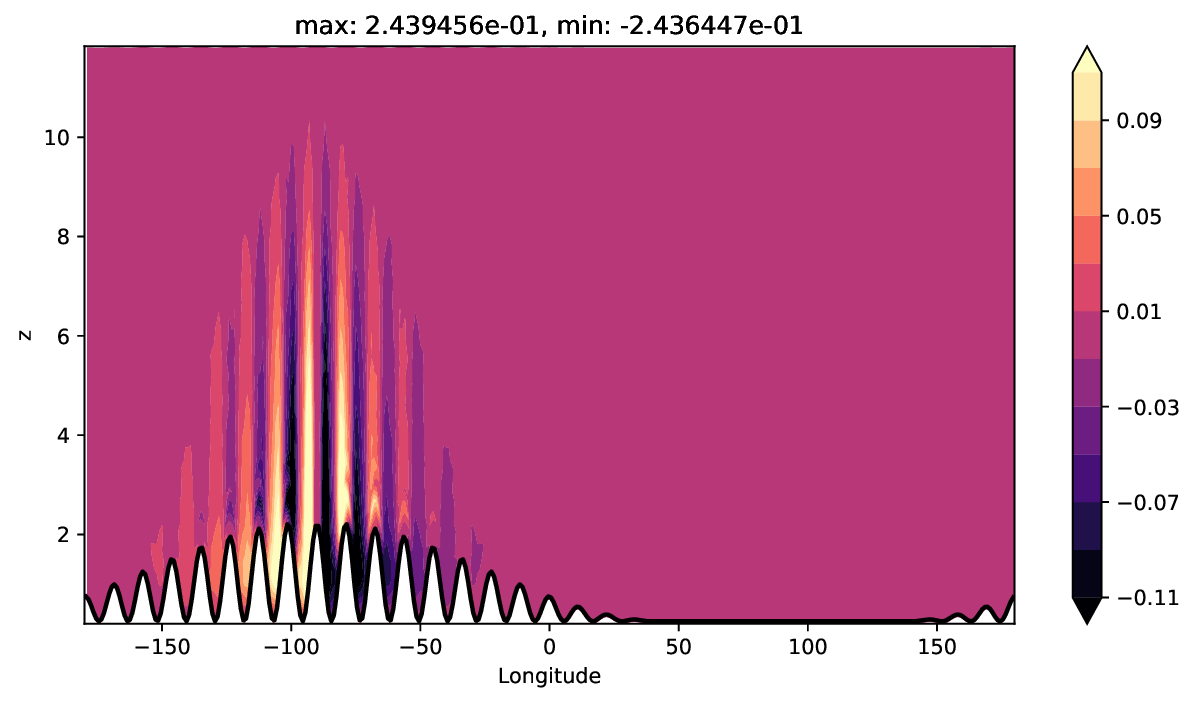}}%
\subfloat[SISL $w$]{\includegraphics[scale=0.3]{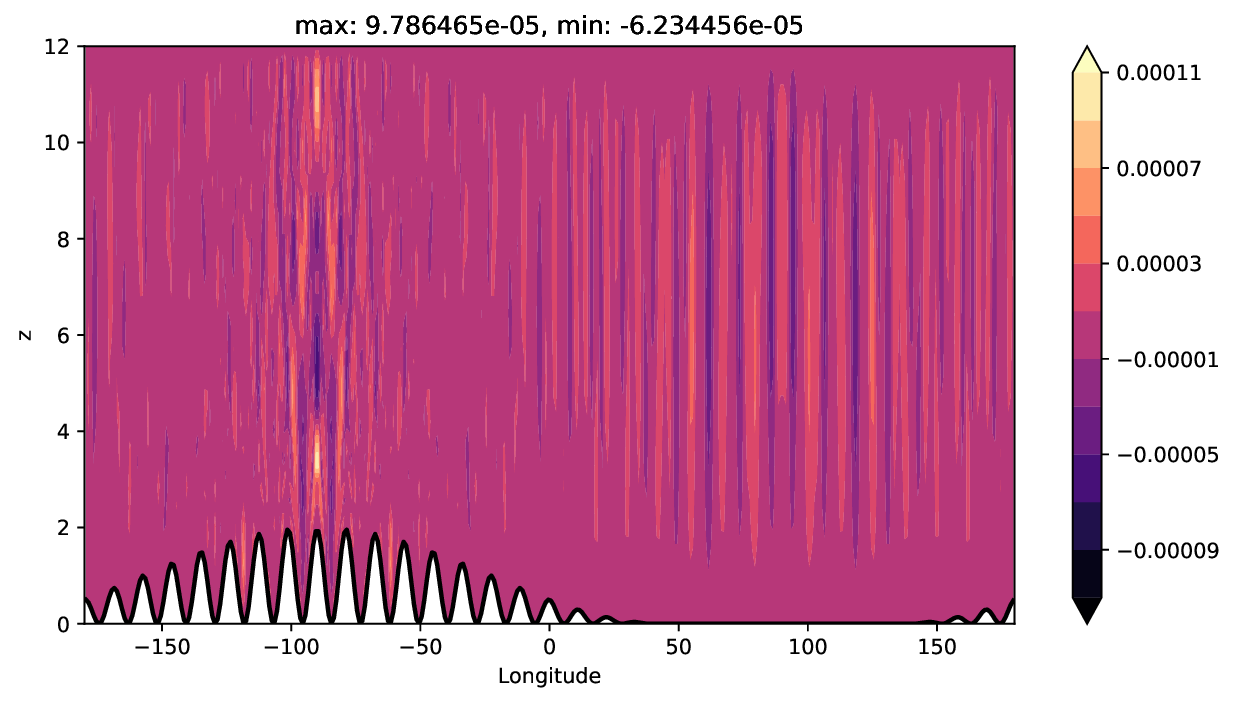}}
\caption{Zonal (left panels) and vertical (right panels) wind fields for the resting atmosphere
case after 6 days on a C96L30 mesh with $\Delta t = 600$ s. The top row shows results from this paper and the bottom row shows results from the semi-implicit semi-Lagrangian ENDGame model with a 1 degree resolution and $\Delta t = 600$ s, (\citet{wood14}).}
\label{fig:dcmip200}
\end{figure}
The test is based upon earlier ideas by \citet{lin97}. The orographic profile is given by
\begin{equation}
    z_B = \begin{cases} 
    \frac{h_0}{2}\left[1+\cos\left(\frac{\pi r_m}{R_m}\right)\right]\cos^2\left(\frac{\pi r_m}{\zeta_m}\right),\quad &\textrm{if}\, r_m < R_m,\\
    0,\quad &\text{otherwise}. 
    \end{cases}\label{eq:dcmip200-zs}
\end{equation}
The mountain height is $h_0 = 2000 m$ and $R_m = 3\pi/4$, $\zeta_m = \pi/16$. The great circle distance from
the mountain centrepoint $\left(\lambda_m,\, \phi_m\right) = \left(3\pi/2, \, 0\right)$ is
\begin{equation}
    r_m = \arccos\left[\sin\phi_m\sin\phi + \cos\phi_m\cos_\phi\cos\left(\lambda-\lambda_m\right)\right].\label{eq:dcmip200-rm}
\end{equation}
The atmosphere is initialised at rest ($\mathbf{u}=0$) and is non-rotating ($\boldsymbol{\Omega}=0$). A constant lapse
rate $\Gamma = 0.0065 K m^{-1}$ is used, giving the initial temperature as 
\begin{equation}
    T = T_0 - \Gamma \left(r-a\right),\label{eq:dcmip200-T}
\end{equation}
with $T_0 = 300 K$. This example tests the accuracy of the pressure gradient terms over orography. Since the initial state is in balance no motion should be generated, but, due to inaccuracies in the pressure gradient terms in terrain following coordinates, this balance will not be discretely maintained and motion will be generated. The size of the motion generated is related to the error in the pressure gradient terms. Figure \ref{fig:dcmip200} shows the zonal and vertical wind fields along the equator
after 6 days simulation on a C96L30 mesh with a uniform vertical resolution and model top at $z_{T} = 12 Km$. This is compared to the semi-implicit semi-Lagrangian model of \citet{wood14} run on a 1 degree Latitude-Longitude mesh with the same vertical mesh.
Over the mountain a small amount of motion is generated, as shown in both the zonal and vertical velocities (Figure \ref{fig:dcmip200}, top row). Compared to \citet{wood14} (Figure \ref{fig:dcmip200}, bottom row) the zonal velocity perturbations are approximately an order of magnitude smaller while the vertical velocity perturbations are of the same order but less widespread, particularly over the orography. Taken together this indicates the model is able to maintain balance over orography relatively well and there is no spuriously large growth of perturbations.


\subsection{Flow over a Gaussian hill}
\label{sec:gaussian}
\begin{figure}[htbp]
\centering
\subfloat[Day 5 $z700$]{\includegraphics[scale=0.15]{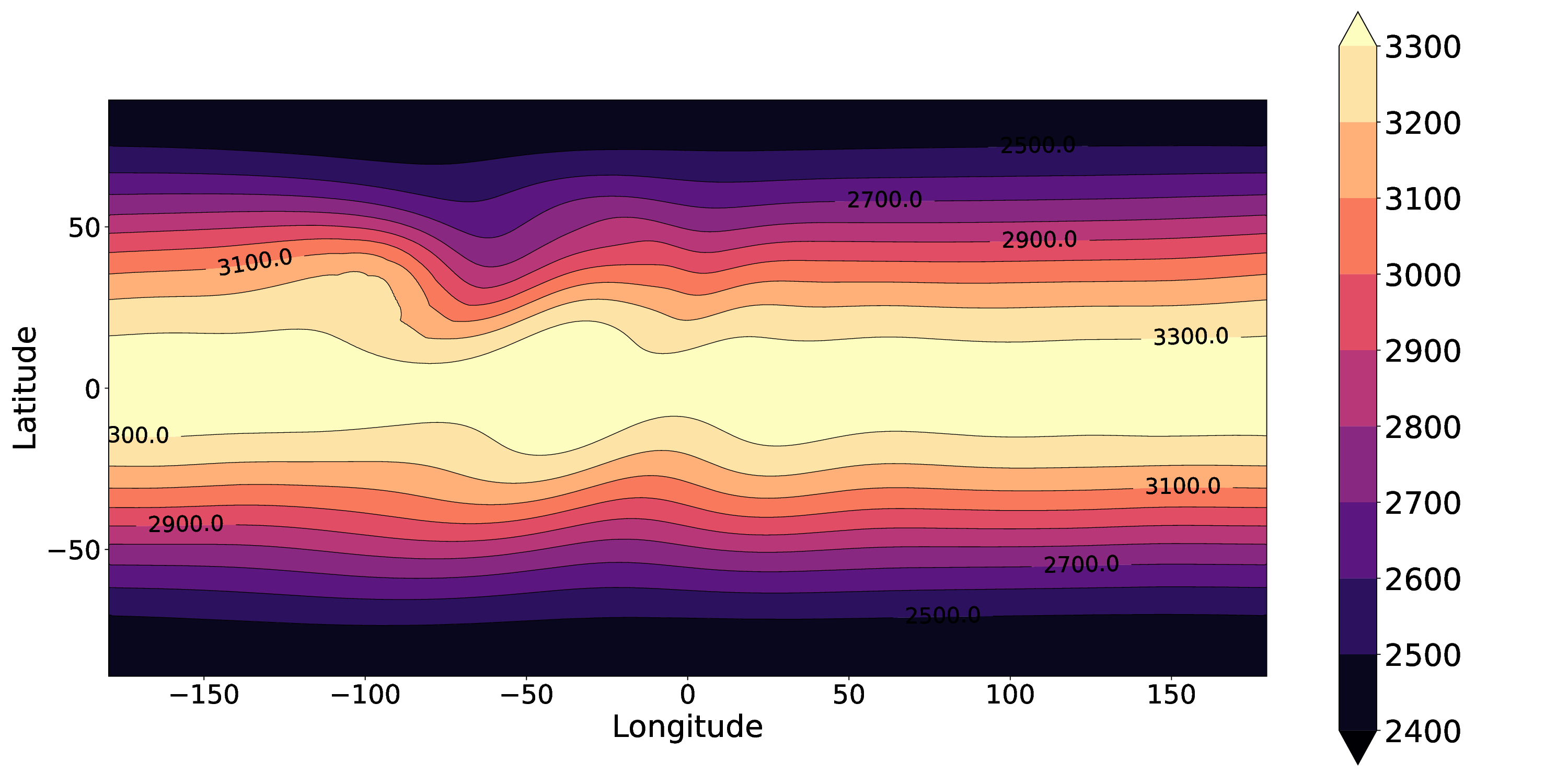}}
\subfloat[Day 5 $T700$]{\includegraphics[scale=0.15]{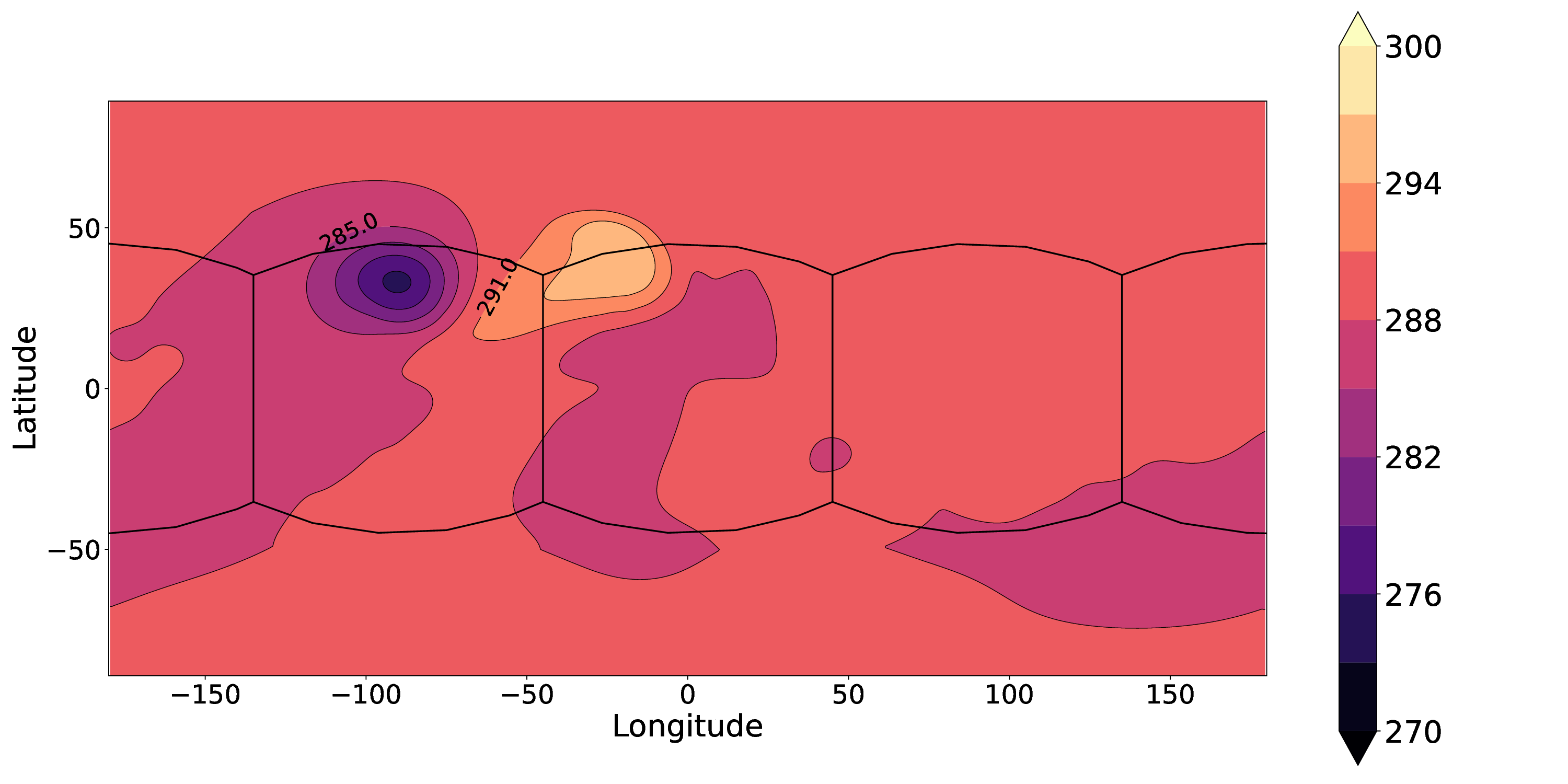}}\\
\subfloat[Day 10 $z700$]{\includegraphics[scale=0.15]{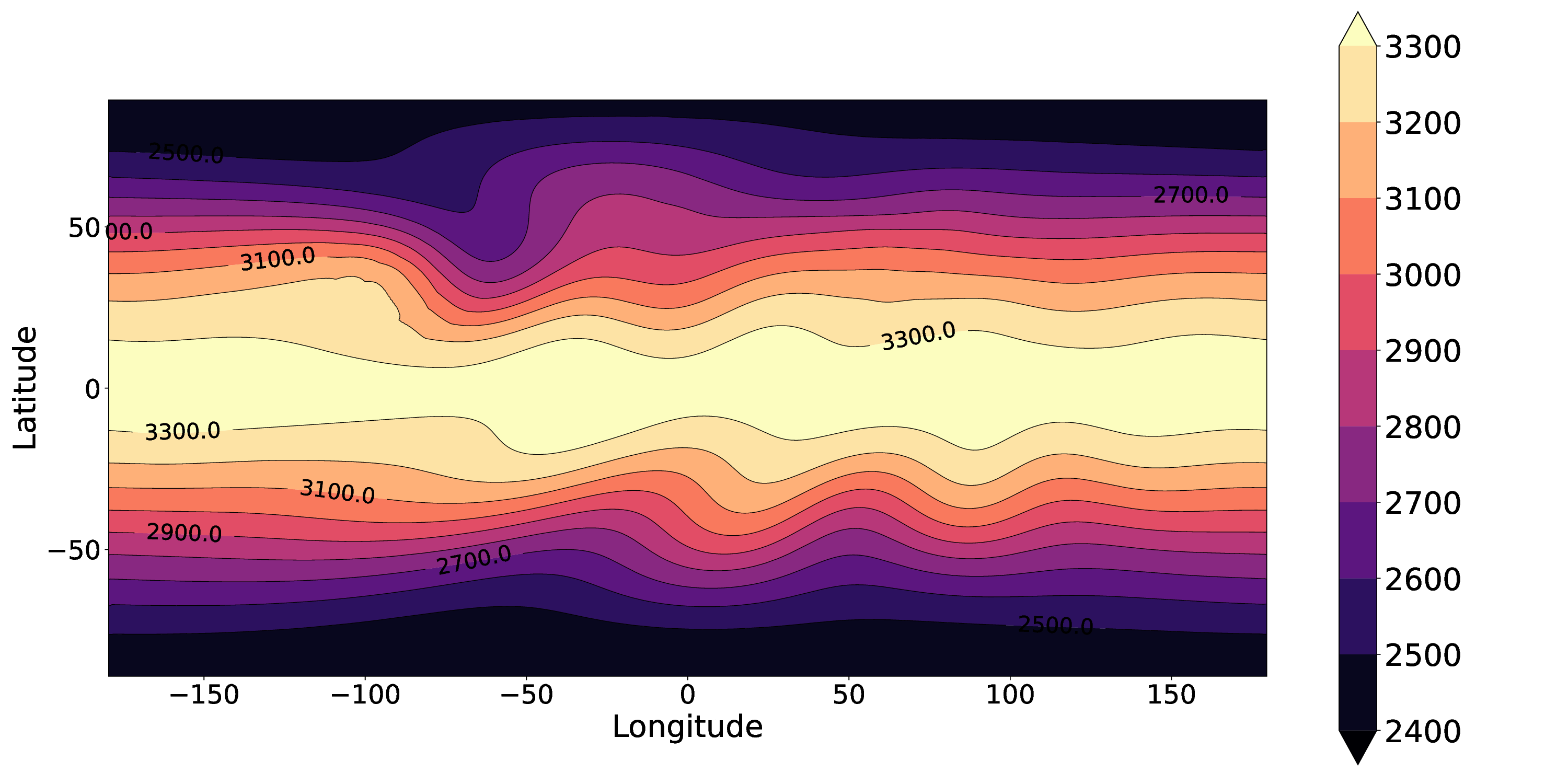}}
\subfloat[Day 10 $T700$]{\includegraphics[scale=0.15]{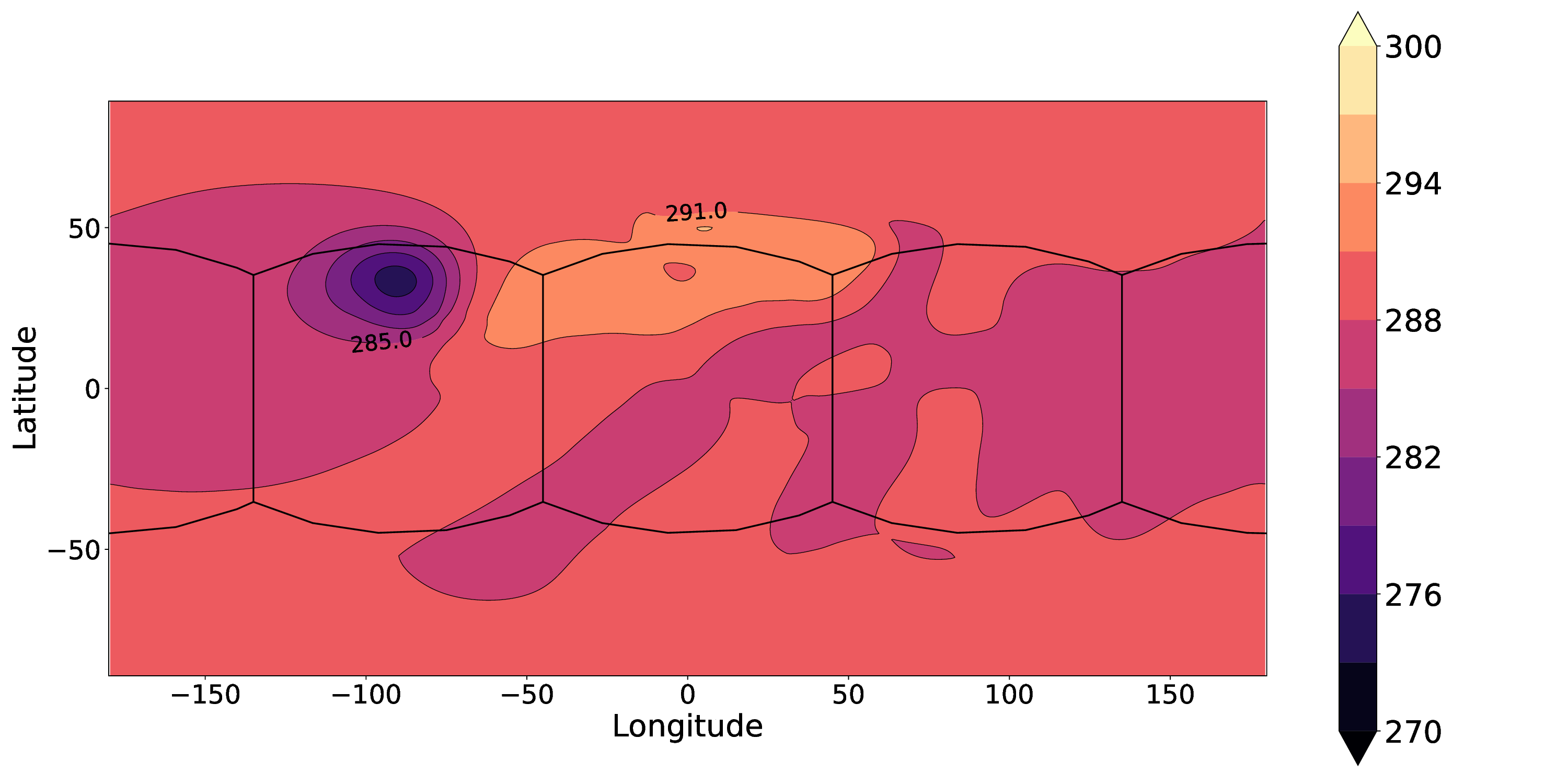}}\\
\subfloat[Day 15 $z700$]{\includegraphics[scale=0.15]{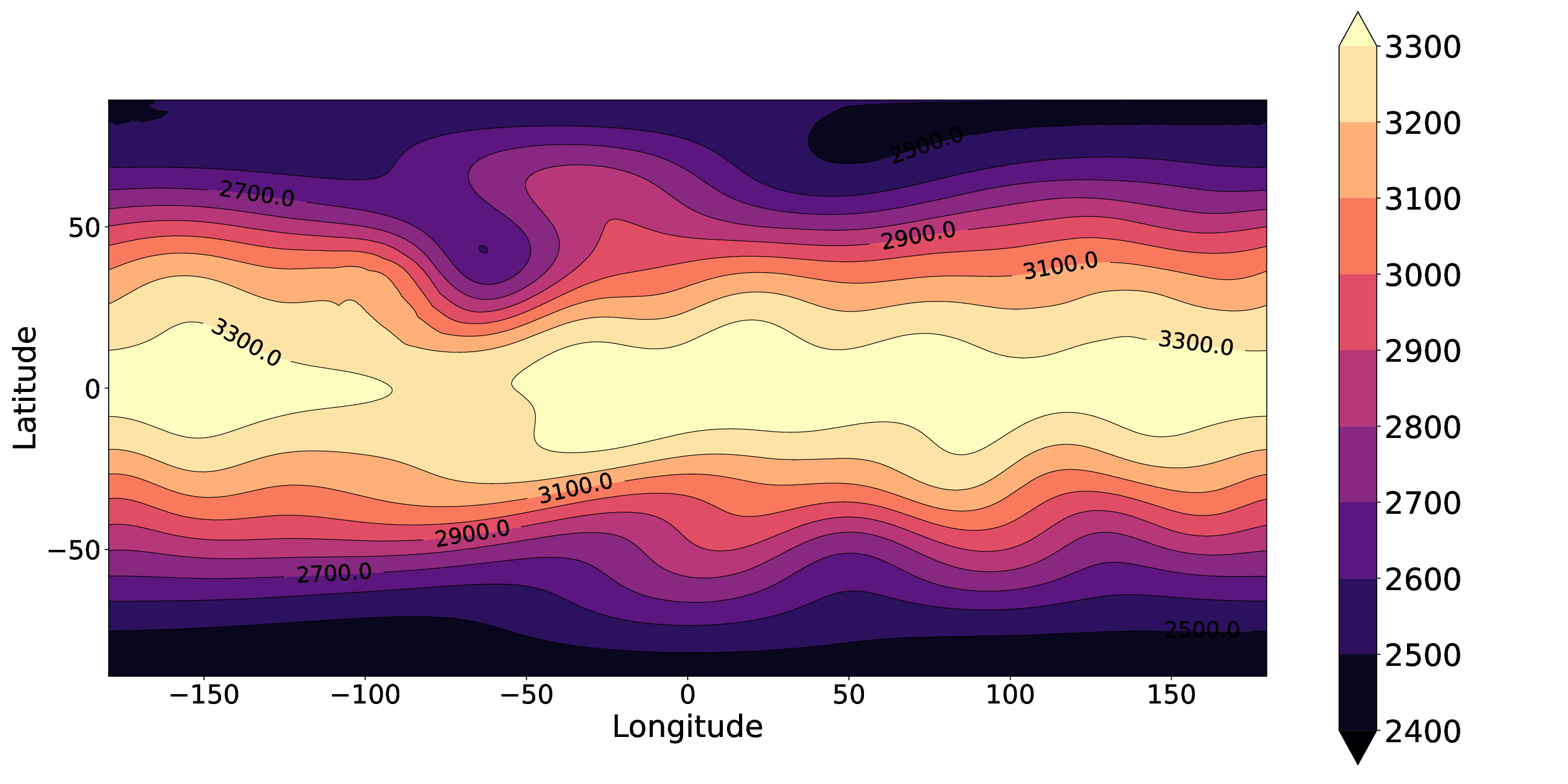}}
\subfloat[Day 15 $T700$]{\includegraphics[scale=0.15]{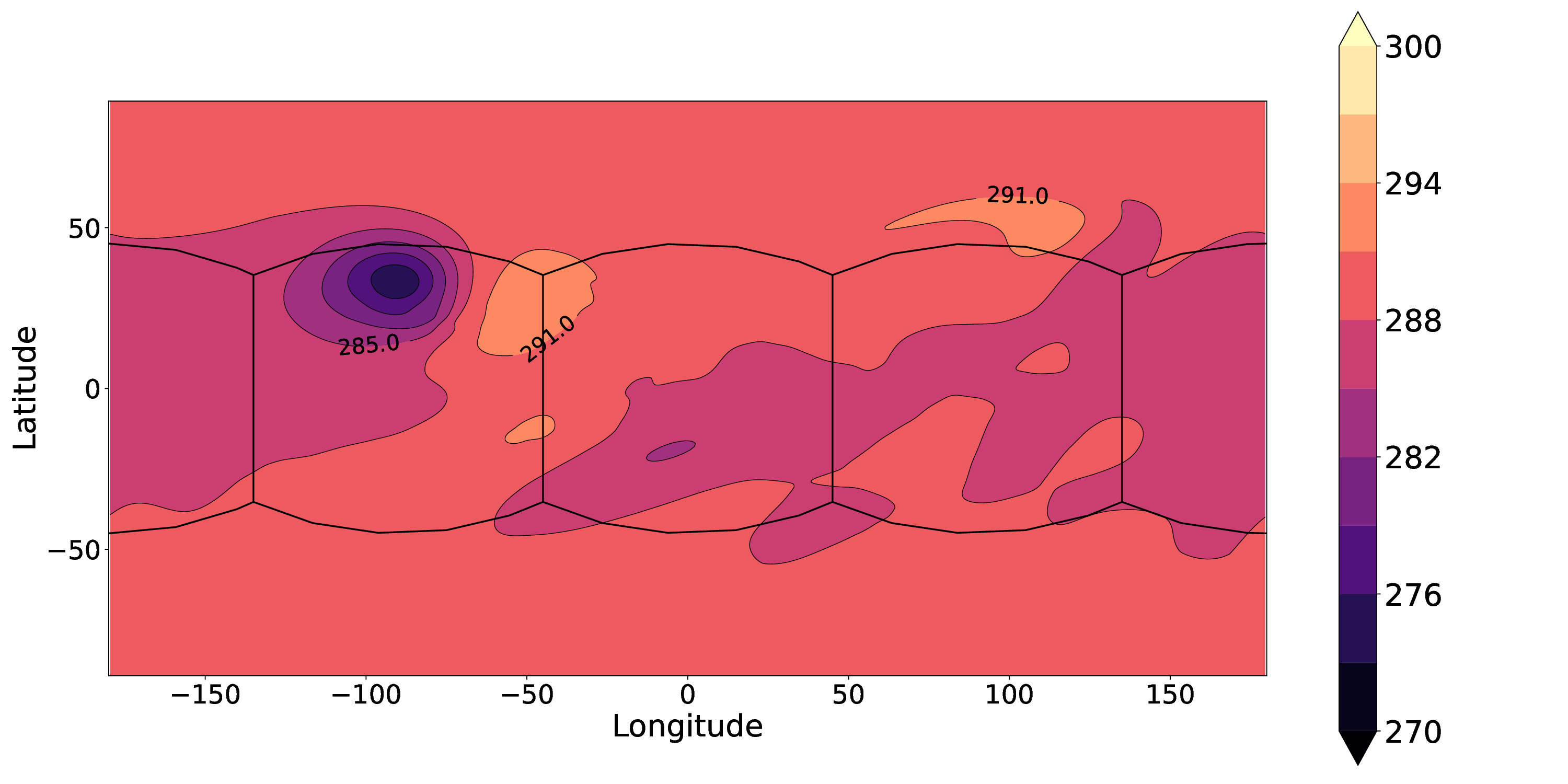}}\\
\caption{$700 hPa$ geopotential height (left column) and temperature (right column) at days 5, 10 and 15 for the flow over a Gaussian hill test at C96L40 resolution with $\Delta t = 900 s$. The location of the cubed-sphere panel boundaries are overlaid on the $T700$ plots.}
\label{fig:sbr-gaussian}
\end{figure}
This test simulates the  generation of Rossby waves from flow over orography and is based upon the test of \citet{tomita04} and \citet{jablonowski08} and further developed by \citet{allen16}.
The set up used here follows \citet{allen16} and the mountain profile is given by
\begin{equation}
    z_B = h_0\exp\left[-\left(\frac{a}{\zeta_m}r_m\right)^2\right],\label{eq:gaussian-h0}
\end{equation}
with $h_0 = 2000 m$, $\zeta_m = 1500 km$ and $r_m$ given by \eqref{eq:dcmip200-rm} with the mountain centre at $\left(\lambda_m,\, \phi_m\right) = \left(-\pi/2, \pi/6\right)$. The initial wind is given by
\begin{equation}
u\left(r\cos\phi\right)=u_0\frac{r}{a}\cos\phi,
\label{eq:whitestan-sbr-wind}
\end{equation}
with $u_0 = 20 ms^{-1}$ which is (5.7) from \citet{allen16} with $\beta=0$. The atmosphere is isothermal with $T = 288 K$ and the surface pressure is given by (5.15) in \citet{allen16}:
\begin{equation}
    p_s = p_p \exp\left[\left(2\Omega a+u_0\right)\frac{u_0}{2R T}\cos^2\phi\right]\exp\left( -\frac{z_s g}{R T} \right),\label{eq:gaussian-p_s}
\end{equation}
with $p_p = 930 hPa$. The model is run at C96L40 resolution with a uniform vertical mesh and model top at $z_{T} = 32 Km$. The $700 hPa$ geopotential height and temperature are shown at days 5, 10 and 15 of the simulation in Figure \ref{fig:sbr-gaussian}. These results are similar to both the Yin-Yang and ENDGame semi-implicit semi-Lagrangian results presented in \citet{allen16} (their figures 8 and 9). The locations of the cubed-sphere panel boundaries are overlaid on the temperature figures and although the mountain is located across a panel edge there are no obvious indications of grid imprinting from the cubed-sphere on the solution profiles.
%
%

\subsection{Baroclinic wave}
\label{sec:baroclinic}
\begin{figure*}[htbp]
\centering
\subfloat[Day 8 surface pressure]{\includegraphics[scale=0.15]{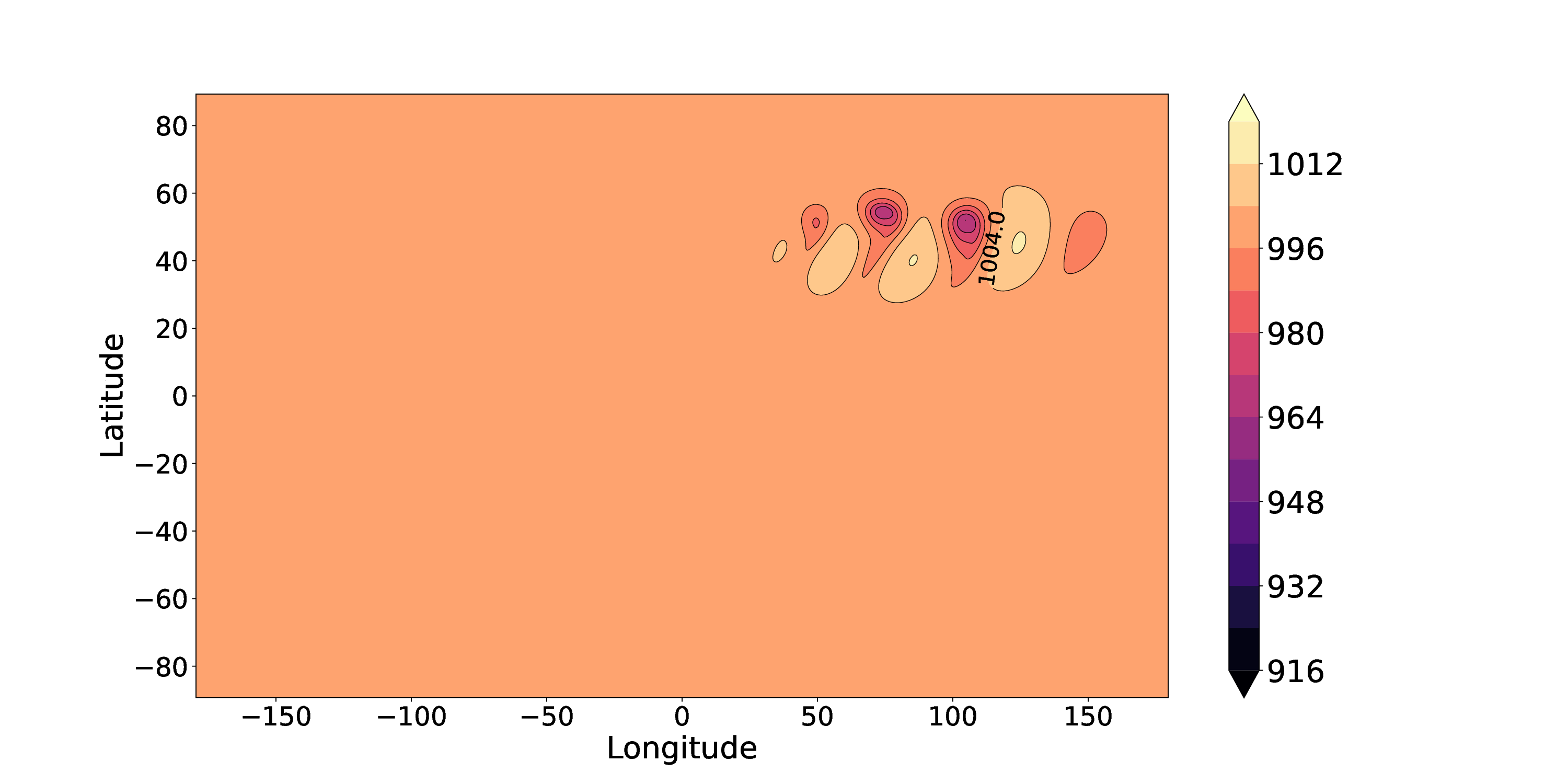}}
\subfloat[Day 10 surface pressure]{\includegraphics[scale=0.15]{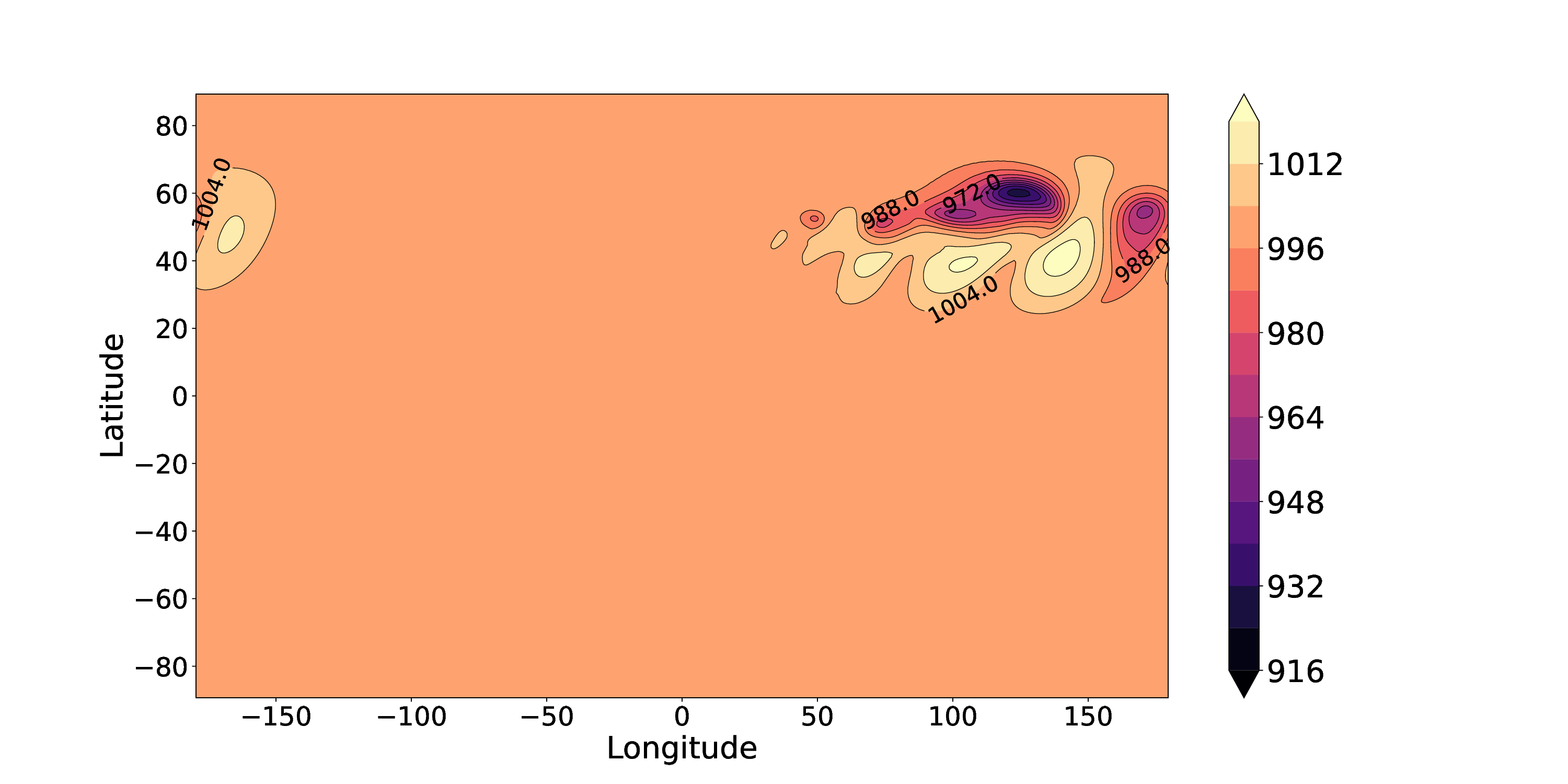}}\\
\subfloat[Day 8 $T850$]{\includegraphics[scale=0.15]{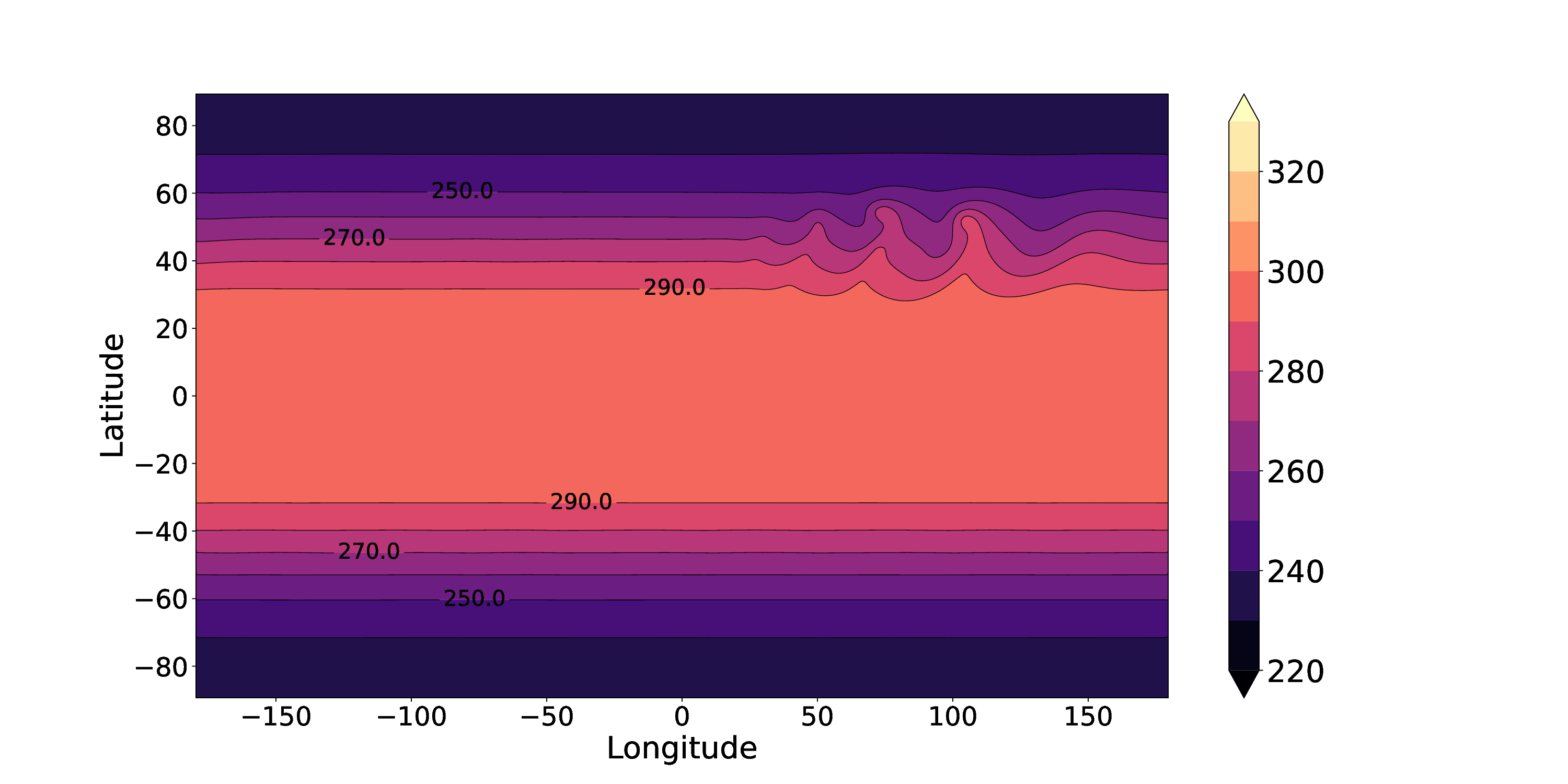}}
\subfloat[Day 10 $T850$]{\includegraphics[scale=0.15]{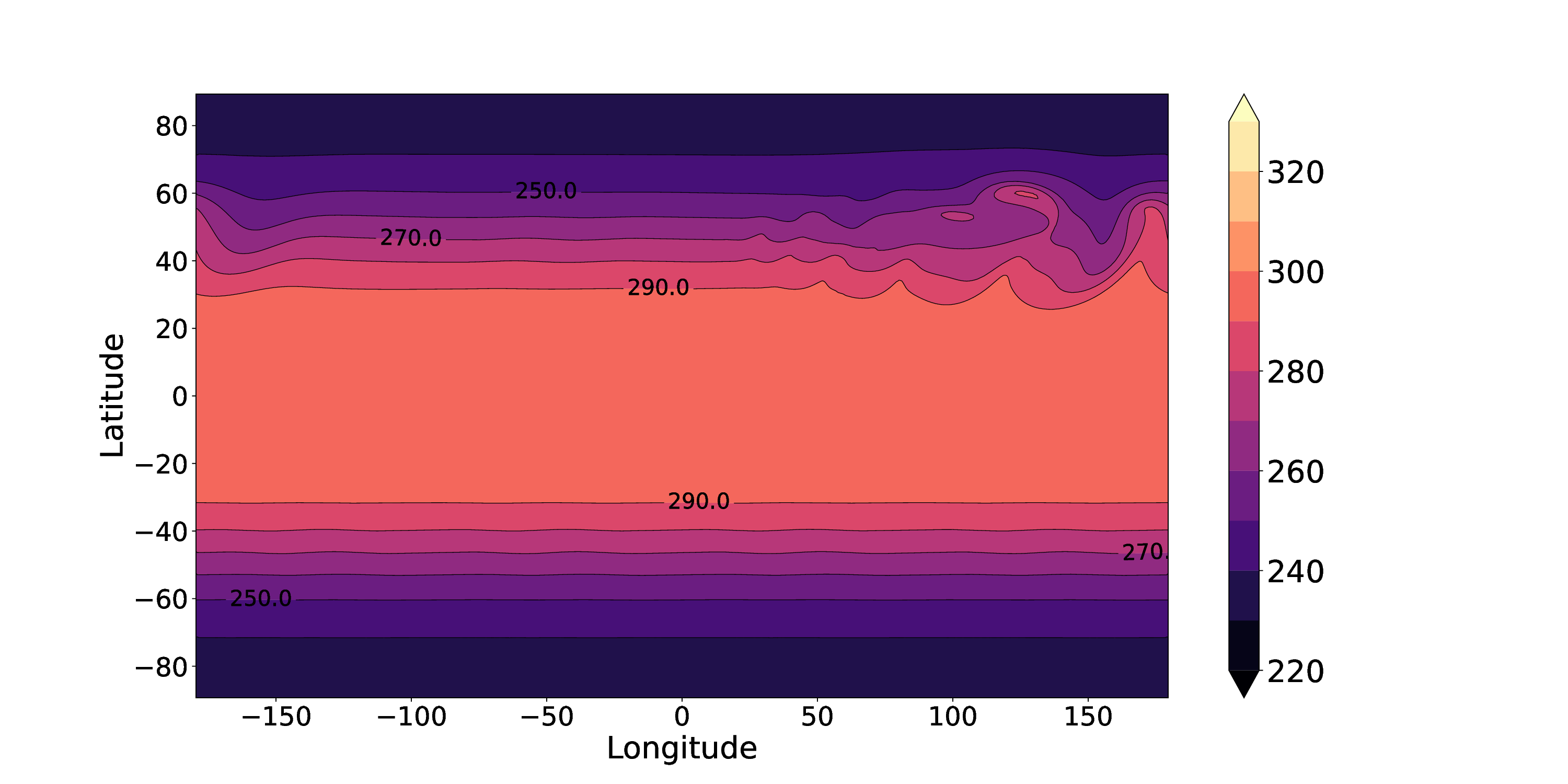}}
\caption{Surface pressure (top row) and 850 hPa temperature (bottom row) for the Baroclinic wave test on a C96L30 mesh with $\Delta t = 900$ s. Left panels: after 8 days simulation and Right panels: after 10 days simulation.}
\label{fig:baroclinic}
\end{figure*}
\begin{figure*}[htbp]
\centering
\includegraphics[scale=0.25]{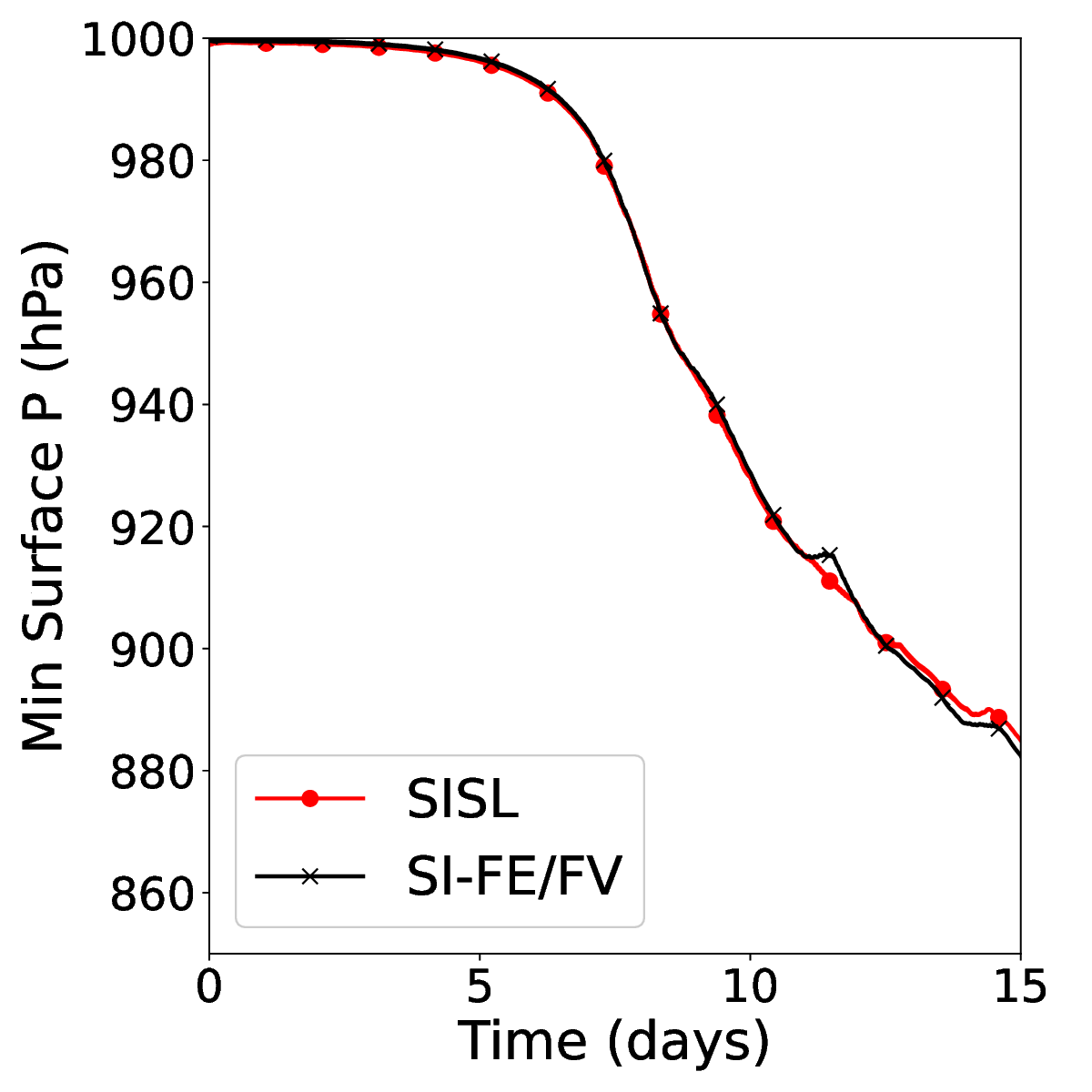}%
\caption{Minimum surface pressure for the Baroclinic wave test on a C96L30 mesh with $\Delta t = 900$ s compared with the SISL model of \citet{wood14}  on a 1 degree mesh with the same timestep. }
\label{fig:baroclinic-min-surface-p}
\end{figure*}
The baroclinic wave test of \citet{ullrich14} simulates the formation of a series of features typical of mid-latitude weather systems. The test is run for 15 days at C96L30 resolution with a quadratic vertical mesh and the same timestep $\Delta t = 900 s$ as used in \citet{wood14} with a monotonic filter applied to the transport of $\theta$ as described in Section \ref{sec:adv_inc_comp}. The surface pressure and $850 hPa$ temperature at days 8 and 10 are shown in Figure \ref{fig:baroclinic}. These compare well with the results shown in both \citet{ullrich14} and \citet{wood14}. Importantly there are no obvious signs of grid imprinting in the southern hemisphere in either the pressure or temperature fields.
The minimum surface pressure throughout the simulation compared to \citet{wood14} is shown in Figure \ref{fig:baroclinic-min-surface-p} and both models show excellent agreement through the first 10 days of simulation before some divergence in model solutions over the last 5 days of simulation.
\begin{figure}[htbp]
\centering
\subfloat[C448]{\includegraphics[scale=0.15]{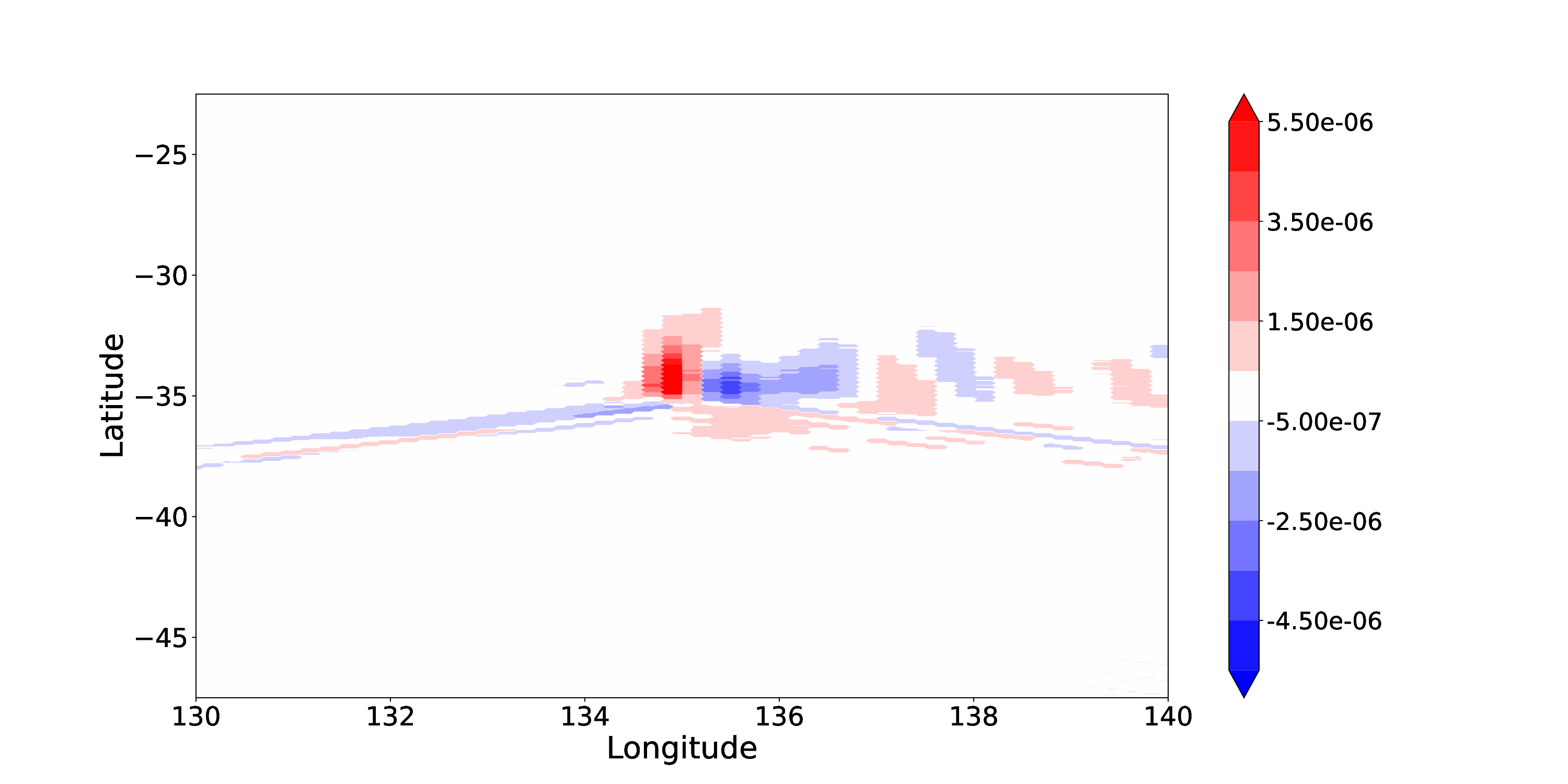}}%
\subfloat[C896]{\includegraphics[scale=0.15]{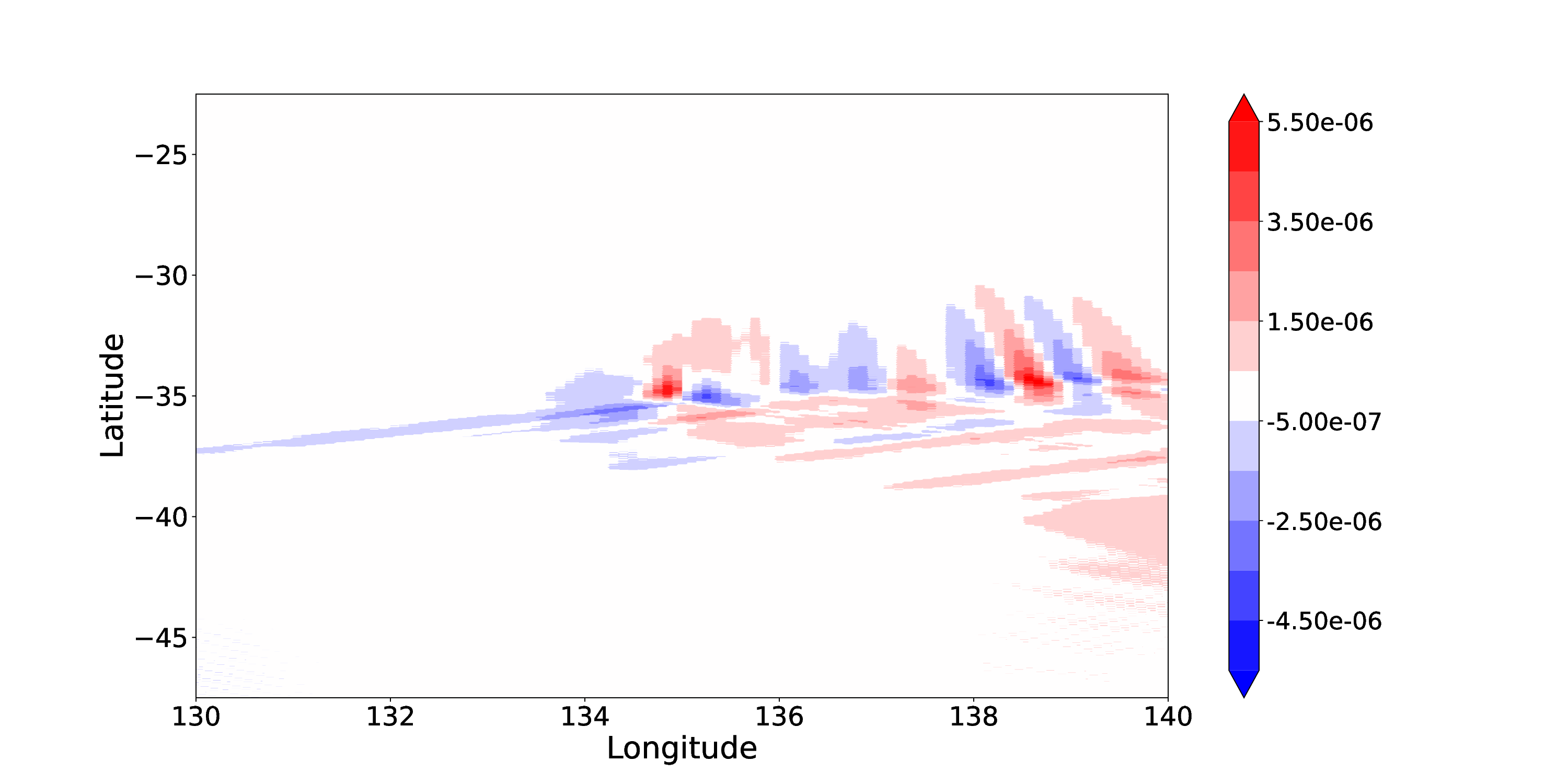}}\\
\caption{Vertical velocity at a corner of the cubed-sphere after 1 day of simulation on a (left) C448L30 mesh with $\Delta t = 450$ and (right) C896L30 mesh with $\Delta t = 225$.}\label{fig:hiwpp}
\end{figure}
\subsubsection{Grid Imprinting}\label{sec:hiwpp}

Following \citet{hiwpp} to further investigate the grid imprinting in the model this test is repeated at high horizontal resolution (C448 and C896) for a single simulated day. The vertical velocity (with zonal mean removed) around a corner of the cubed-sphere in the southern hemisphere (away from the initial perturbation) is shown in Figure \ref{fig:hiwpp} and can be compared with Figures 2.6 \& 2.7 of \citet{hiwpp}. At both resolutions there is some spurious motion around the corner of the cubed-sphere, however the values are very small, of the order $10^{-3}$mm/s, and smaller than both the MPAS and FV3 models considered in \citet{hiwpp}, so this level of grid imprinting is considered acceptable. These results support the choice to use a lowest degree compatible finite-element method  which simplifies a number of aspects of the model design (such as coupling to existing subgrid parametrisation schemes, \citet{brown23}).


\subsection{Held-Suarez}
\label{sec:held-suarez}
\begin{figure*}[htbp]
\centering
\includegraphics[scale=0.13]{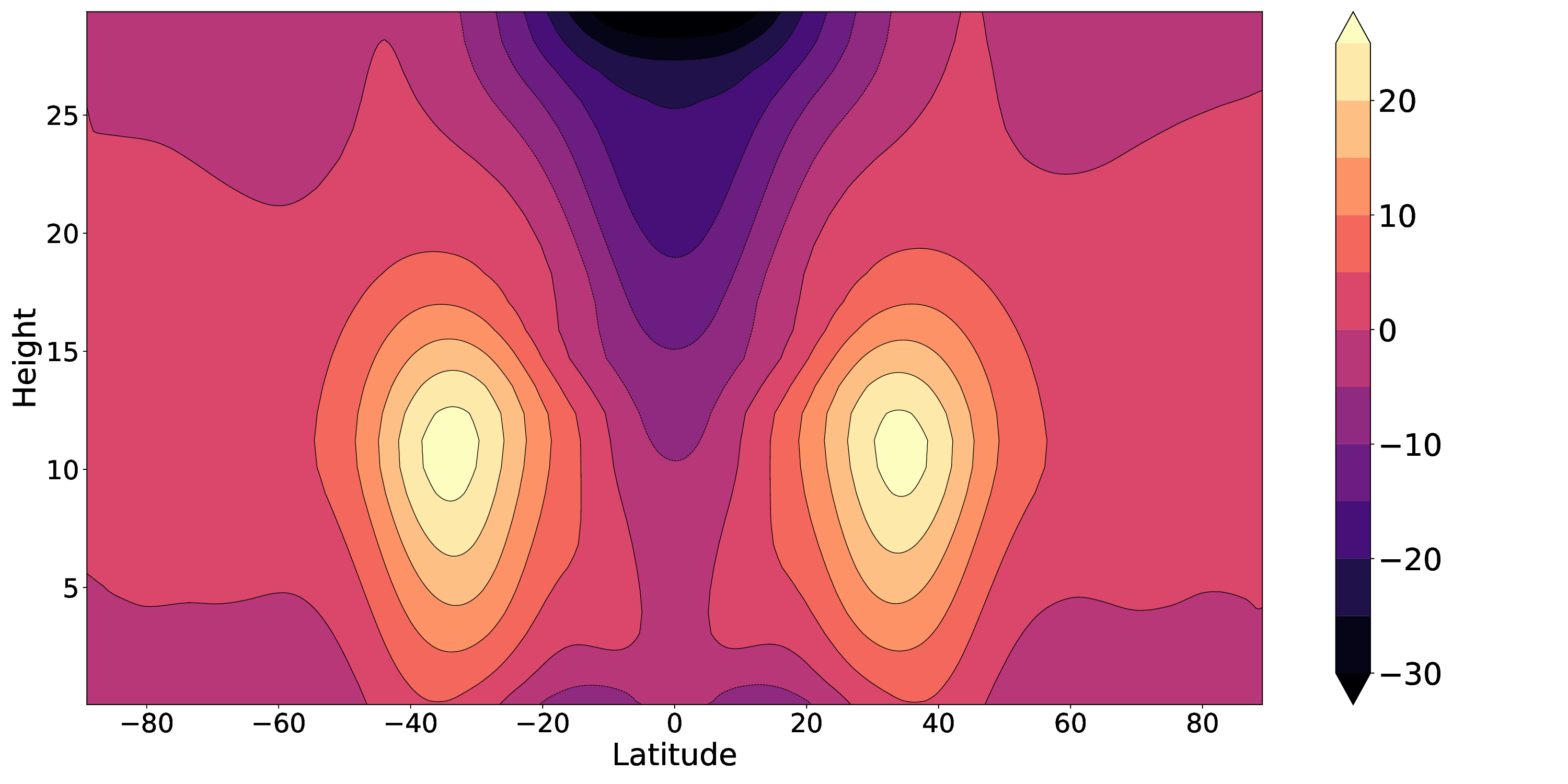}%
\includegraphics[scale=0.13]{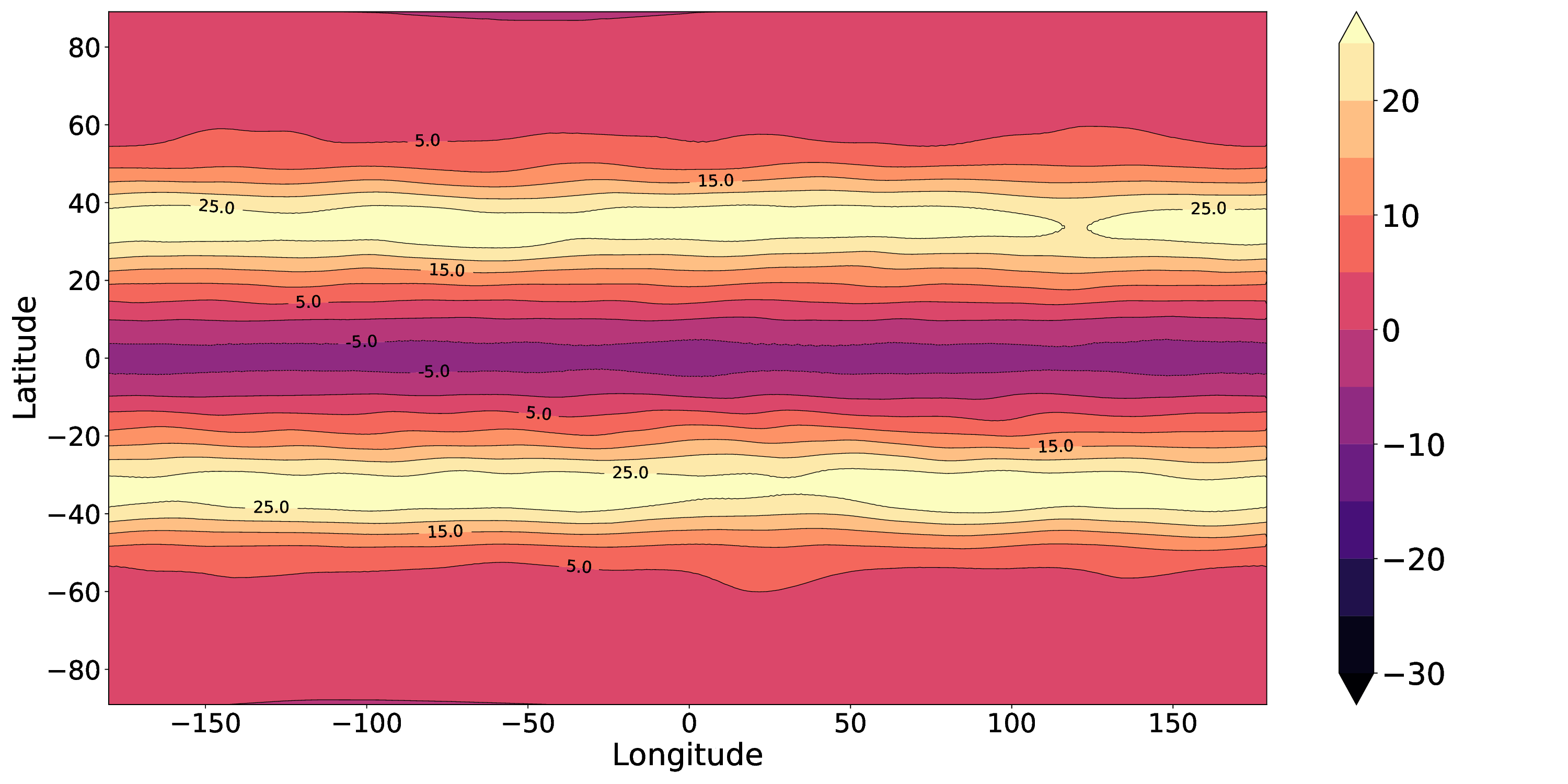}\\
\caption{Zonally averaged zonal velocity (left) and zonal velocity on level 14 ($\approx 10.5 km$) (right) fields. All fields are averaged over the last 1000 days of a 1200 day Held-Suarez run using a C48L30 mesh with $\Delta t = 1800$ s.}
\label{fig:held-suarez}
\end{figure*}

The Held-Suarez climate test (\citet{HeldSuarez94}) simulates the evolution of an atmospheric state with idealised surface temperature and wind forcing. The model is run for 1200 days and time-averaged fields (after the first 200 day spinup period) are shown in Figure \ref{fig:held-suarez}.  The initial state is taken to be the baroclinic wave initial state from Section \ref{sec:baroclinic}.  
The time-average of the zonal velocity
field on a C48L30 mesh using a quadratic stretching in the vertical mesh is shown in Figure \ref{fig:held-suarez}. The left panel shows the zonally averaged zonal velocity field and the right panel shows the zonal velocity on level 14 (approximately though the centre of the jets). These profiles are again similar to those produced by ENDGame (\citet{Tort15}). The horizontal cross section shows no obvious sign of grid imprinting from the cubed-sphere grid, indicating that even for long timescale runs there is no systematic error from the treatment of the cubed-sphere panel corners and edges. 
%

%
%


\section{Summary}
\label{sec:Discussion}

The mixed finite-element, finite-volume, semi-implicit model of \citet{melvin19} has been extended to spherical geometry for atmospheric modelling. The finite-volume transport scheme has been extended to encompass the momentum advection terms and has been adapted to the non-uniform, non-orthogonal horizontal mesh by using  polynomial reconstructions as in \citet{kent22}. In order to maintain a constant field, an advective-then-flux formulation for conservatively transported variables is used (\citet{bendall23}). The finite-element method presented in \citet{melvin19} has required minimal modification for spherical geometries, highlighting the flexibility of this approach, and the only significant change here is the use of a semi-analytic mapping to the sphere to accurately represent a spherical shell domain. The temporal discretisation mimics that of \citet{wood14} and is coupled to the multigrid solver of \citet{maynard20} to give an efficient method of solving the linear semi-implicit system.

The semi-implicit finite-element finite-volume dynamical core has been applied to a number of standard dynamical tests. The model has been shown to produce results comparable to those in the literature and in particular to the existing semi-implicit semi-Lagrangian ENDGame dynamical core (\citet{wood14}) used in operational models at the Met Office. This model has also recently been coupled to idealised physical parametrisations and chemical processes by \citet{brown23} and to the simulation of Exoplanet atmospheres, \citet{sergeev23}, again producing results comparable to those in the literature.

Continual improvements are being made to optimise the computational performance of the model. There does not appear to be any fundamental barrier to achieving comparable throughput to the present operational model.

As a next step towards using this model for numerical weather and climate prediction the dynamical core has been extended to handle moist dynamic processes and coupled to the existing suite of subgrid physical parametrisations used by \citet{walters19}. The formulation and code have been extended to limited area domains with forced boundaries and both of  these developments will be reported upon in future work. In order to couple with a data assimilation system a tangent linear version of this model has also been developed and again will be reported upon in future work.


\section*{Acknowledgments}
The authors would like to express their gratitude to the many people who have worked and contributed to the development of the LFRic atmosphere model, PSyclone and the GungHo dynamical core, without whose efforts this manuscript would not have been possible.

The work of John Thuburn was supported by NERC under the GungHo Phase 2 programme, grant NE/K006762/1.
The work of Colin Cotter was supported by NERC under the grant NE/K006789/1 and by EPSRC by the grants EP/L016613/1 \& EP/R029423/1.

For the purpose of open access, the author has applied a ‘Creative Commons Attribution' (CC BY) licence to any Author Accepted Manuscript version arising'

\bibliography{main}
\appendix

\end{document}